%% file: ms.tex
\documentclass[10pt]{article}          
\usepackage{graphicx}
\usepackage{url}
\usepackage{amssymb,amsthm}

\usepackage[caption=false]{subfig}
\usepackage{tikz,pgfplots}

\usepackage{bbm}
\usepackage{amsmath,amsfonts,bm} 
\usepackage{algorithm}
\usepackage{algpseudocode}
\usepackage[numbers]{natbib}

\usepackage{geometry}
\def\Argmin{\mathop{\mathrm{Argmin}}}
\def\argmin{\mathop{\mathrm{argmin}}}
\def\diag{\mathrm{diag}}
\def\R{\mathbb{R}}
\def\C{\mathbb{R}}

\def\one{\mathbbm{1}}

\def\bb{\mathbf{b}}

\def\bw{\mathbf{w}}
\def\bx{\mathbf{x}}
\def\by{\mathbf{y}}
\def\bz{\mathbf{z}}
\def\blambda{\bm{\lambda}}
\def\bX{\mathbf{X}}
\def\bY{\mathbf{Y}}
\def\bW{\mathbf{W}}
\def\bL{\mathbf{L}}

\def\bpsi{\bm{\psi}}
\def\bmu{\bm{\mu}}
\newcommand{\red}[1]{\textcolor{red}{#1}}

\definecolor{figure}{rgb}{0.1,0.1,0.6}
\definecolor{myred}{rgb}{0.7,0.05,0.2}
\definecolor{mygreen}{rgb}{0.1,0.7,0.1}
\def\red#1{\textcolor{myred}{#1}}

\def\green#1{\textcolor{mygreen}{#1}} 

\newcommand\fdg[1]{{#1}}

\newtheorem{assumption}{Assumption}
\newtheorem{theorem}{Theorem}
\newtheorem{definition}{Definition}
\newtheorem{proposition}{Proposition}
\newtheorem{remark}{Remark}
\begin{document}

\title{Optimal Transport Approximation of 2-Dimensional Measures}


\author{ Fr\'ed\'eric de Gournay  \and Jonas Kahn \and L\'eo Lebrat \and Pierre Weiss}


%

\maketitle
\begin{abstract}
  We propose a fast and scalable algorithm to project a given density on a set of structured measures  defined over a compact 2D domain. The measures can be discrete or supported on curves for instance. The proposed principle and algorithm are a natural generalization of previous results revolving around the generation of blue-noise point distributions, such as Lloyd's algorithm or more advanced techniques based on power diagrams.
  We analyze the convergence properties and propose new approaches to accelerate the generation of point distributions.
  We also design new algorithms to project curves onto spaces of curves with bounded length and curvature or speed and acceleration.
  We illustrate the algorithm's interest through applications in advanced sampling theory, non-photorealistic rendering and path planning.
\end{abstract}

\begin{figure*}[!h]
\centering
 \subfloat[Original]{\includegraphics[width=0.24\textwidth]{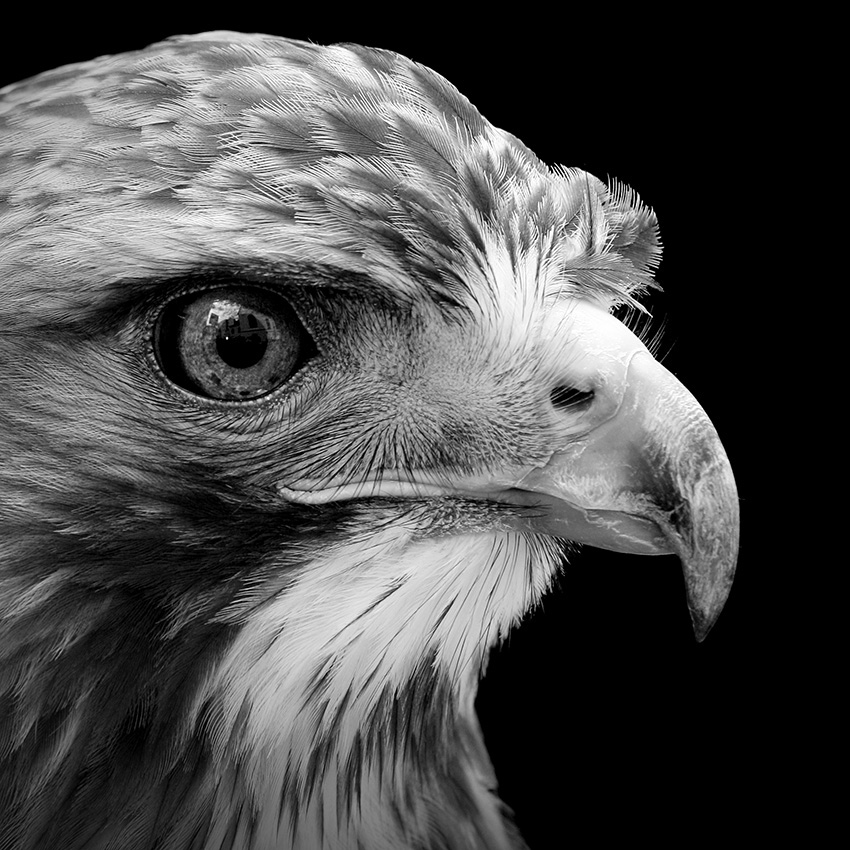}\label{fig:i1}}
 \subfloat[Stippling]{\includegraphics[width=0.24\textwidth]{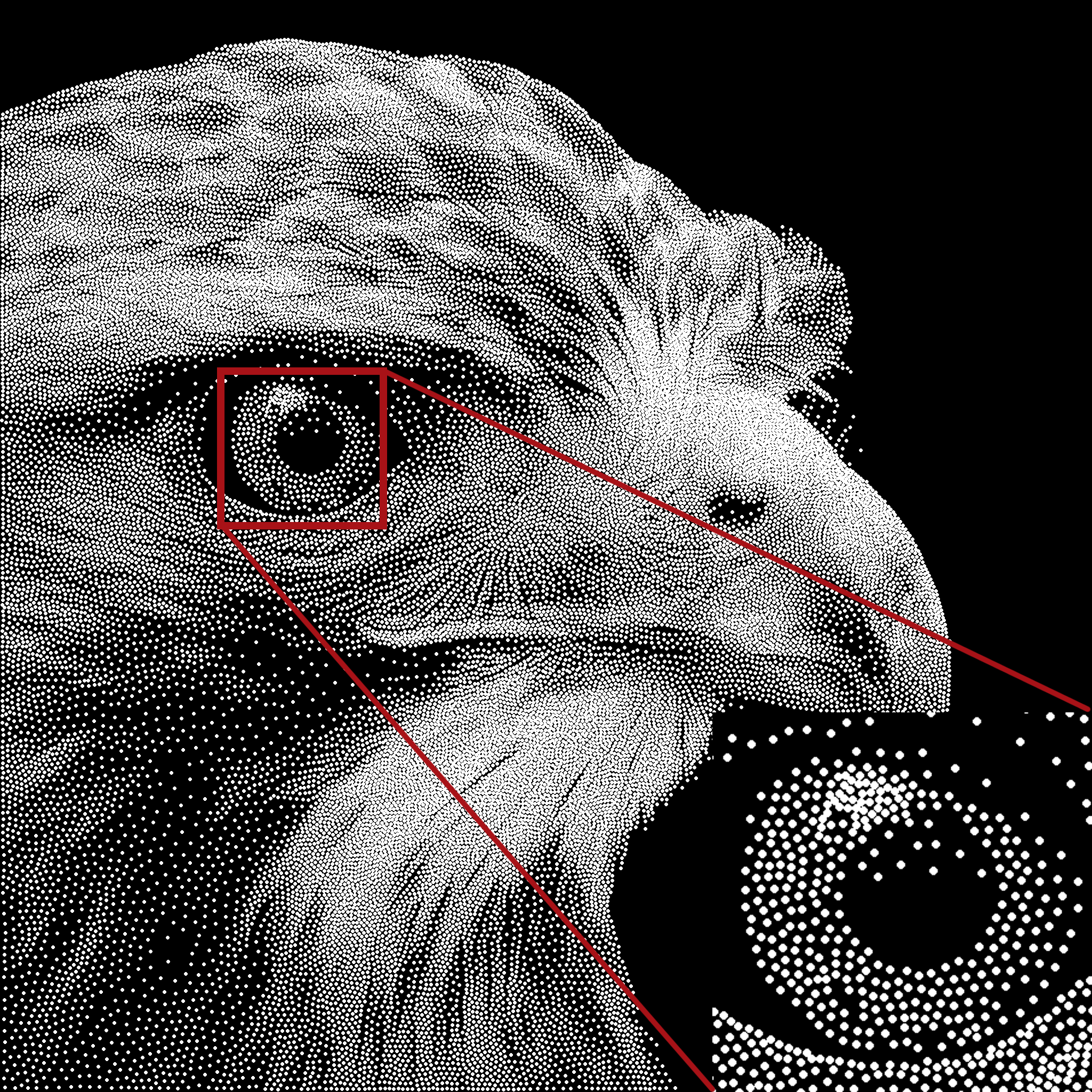}\label{fig:i2}}
 \subfloat[Curvling]{\includegraphics[width=0.24\textwidth]{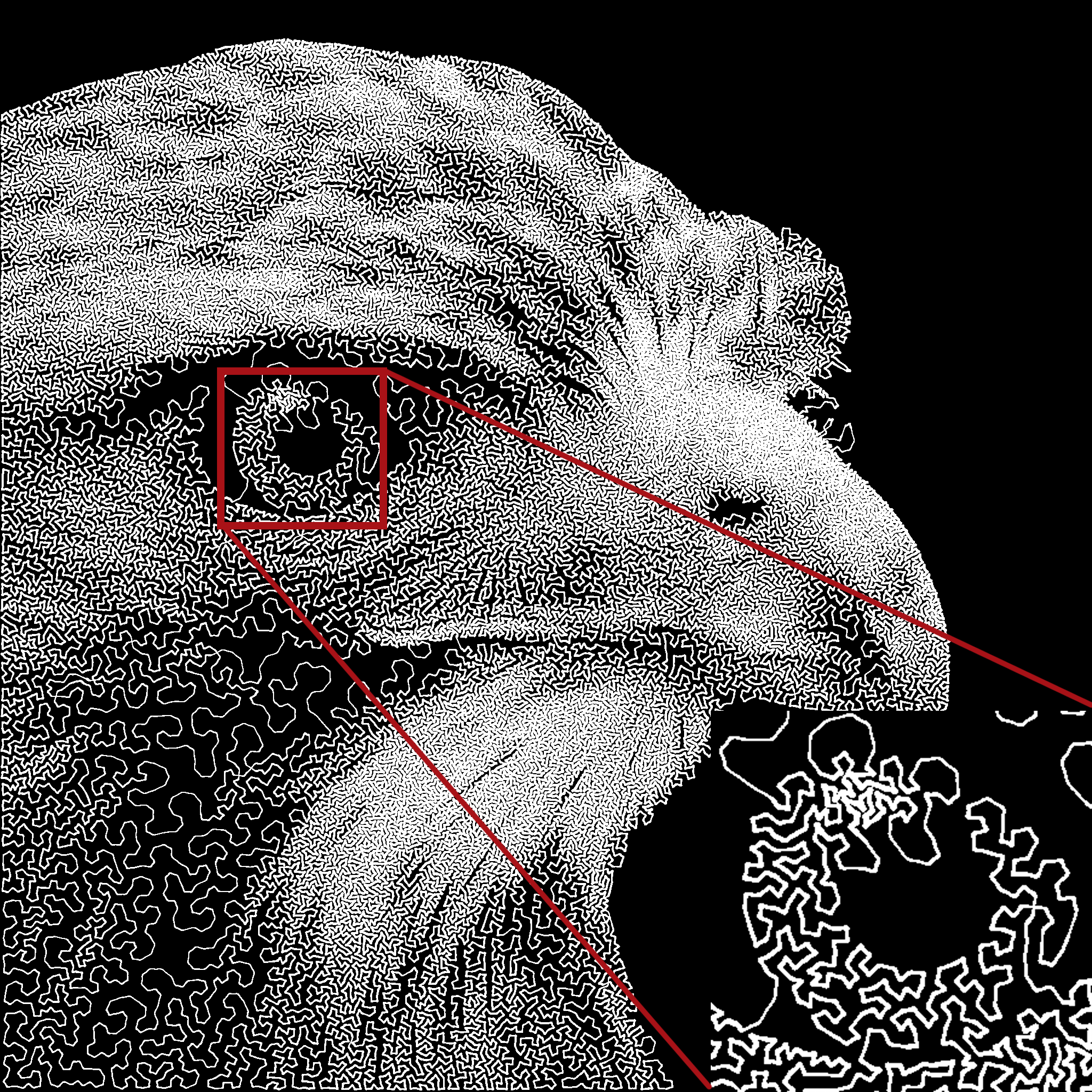}\label{fig:i3}}
 \subfloat[Dashing]{\includegraphics[width=0.24\textwidth]{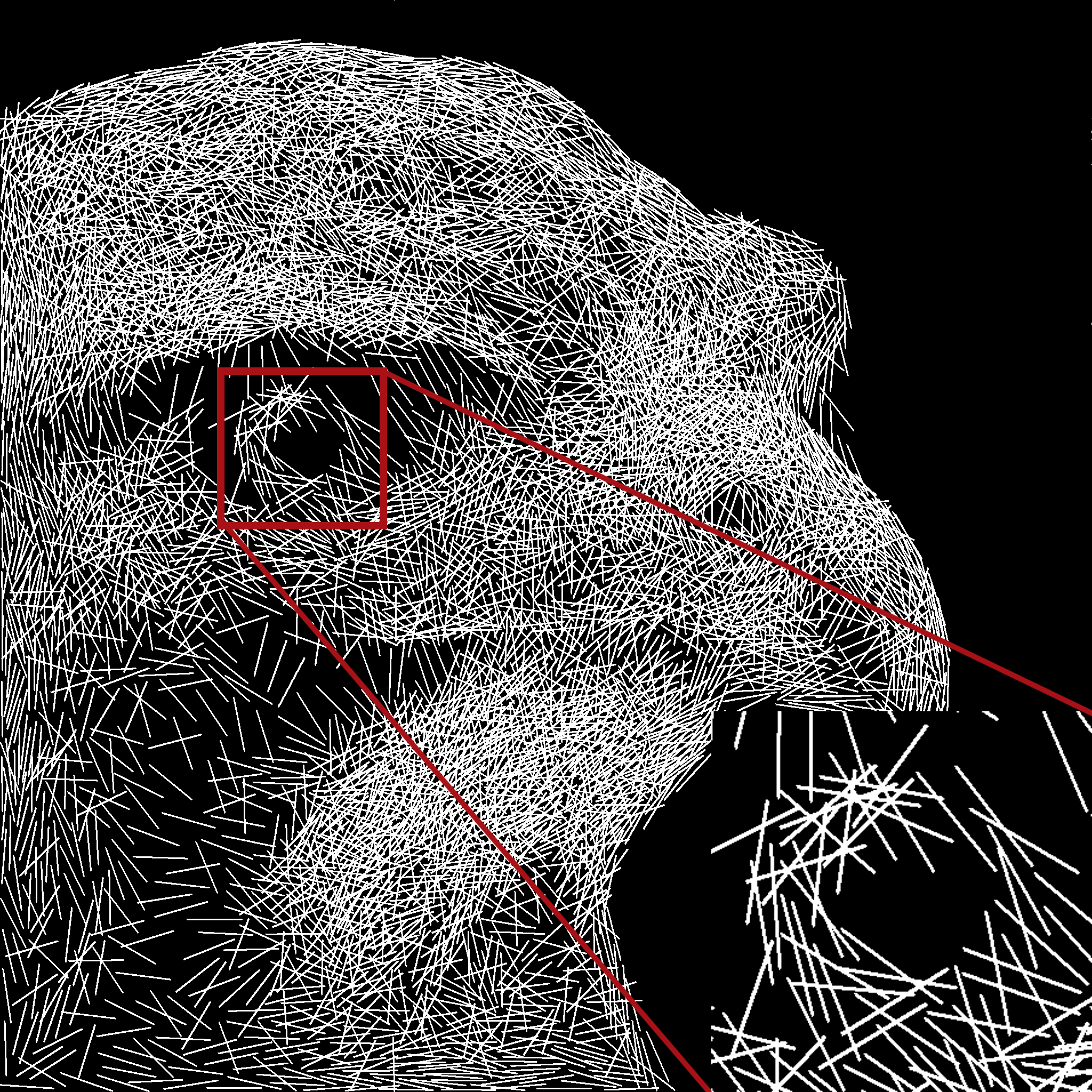}\label{fig:i4}}
 \caption{Approximating an image with a measure supported on points (stippling, 100k, 202''), curve (curvling, 100k, 313'') or segments (dashing, 33k, 237''). In each case, the iterative algorithm starts from a set of points drawn uniformly at random. \label{fig:initial}} 
\end{figure*}

\section{Introduction}

The aim of this paper is to approximate a target measure $\mu$ with probability density function $\rho:\Omega\to\R_+$ with probability measures possessing some structure. This problem arises in a large variety of fields including finance \cite{pages2012optimal}, computer graphics \cite{solomon2015convolutional}, sampling theory \cite{boyer2016generation} or optimal facility location \cite{gastner2006optimal}. 
An example in non photo-realistic rendering is shown in Figure~\ref{fig:initial}, where the target image in Fig.~\ref{fig:i1} is approximated by an atomic measure in Fig.~\ref{fig:i2}, by a smooth curve in Fig.~\ref{fig:i3} and by a set of segments in Fig.~\ref{fig:i4}. 
Given a set of admissible measures $\mathcal{M}$ (i.e. atomic measures, measures supported on smooth curves or segments), the best approximation problem can be expressed as follows:
\begin{equation}\label{eq:mainproblem0}
 \min_{\nu \in \mathcal{M}} D(\nu,\mu),
\end{equation}
where $D$ is a distance between measures.

\subsection{Contributions}

Our main contributions in this article are listed below.
\begin{itemize}
 \item Develop a few original applications for the proposed algorithm.
 \item Develop a fast numerical algorithm to minimize problem \eqref{eq:mainproblem0}, when $D$ is the $W_2$ transportation distance and $\Omega=[0,1]^2$.
 \item Show its connections to existing methods such as Lloyd's algorithm \cite{lloyd1982least} or optimal transport halftoning \cite{DeGoes}.
 \item  Provide some theoretical convergence guarantees for the computation of the optimal semi-discrete transportation plan, especially for complicated point configurations and densities, for which an analysis was still lacking.
 \item Design algorithms specific to the case where the space of admissible measures $\mathcal{M}$ consists of measures supported on curves with geometric constraints (e.g. fixed length and bounded curvature).
 \item Generate a gallery of results to show the versatility of the approach.
\end{itemize}
In the next section, we put our main contributions in perspective.

\subsection{Related works}

\subsubsection{Projections on measure spaces}

To the best of our knowledge, the generic problem \eqref{eq:mainproblem0} was first proposed in \cite{chauffert2017projection} with a distance $D$ constructed through a convolution kernel. 
Similar problems were considered earlier, with spaces of measures restricted to a fixed support for instance \cite{mcasey1999optimal}, but not with the same level of generality.

Formulation \eqref{eq:mainproblem0} covers a large amount of applications that are often not formulated explicitly as optimization problems. We review a few of them below.

\paragraph{Finitely supported measures}

A lot of approaches have been developed when $\mathcal{M}$ is the set of  uniform finitely supported measures
\begin{equation}\label{eq:npoint_measures}
 \mathcal{M}_{f,n}=\left\{\nu(\bx)=\frac{1}{n}\sum_{i=1}^n  \delta_{\bx[i]}, \bx \in \Omega^n\right\},
\end{equation}
where $n$ is the support cardinality, or the set of atomic measures defined by:
\begin{equation}\label{eq:atomic_measures}
	 \mathcal{M}_{a,n}\!\!=\!\left\{\nu(\bx,\bw)\! = \!\sum_{i=1}^n \!\! \bw[i] \delta_{\bx[i]}, \bx \in \Omega^n, \bw \in\Delta_{n-1}\right\}\!,
\end{equation}
where $\Delta_{n-1}=\left\{\sum_{i=1}^n \bw[i] = 1, \bw[i]\geq 0, \ \forall i\right\}$ is the canonical simplex.

For these finitely supported measure sets, solving problem \eqref{eq:mainproblem0} yields nice stippling results, which is the process of approximating an image with a finite set of dots (see Fig.~\ref{fig:i2}). 
This problem has a long history and a large amount of methods were designed to find dots locations and radii that minimize  visual artifacts due to discretization \cite{floyd1976adaptive,lloyd1982least,ulichney1987digital,balzer2009capacity}.
Lloyd's algorithm is among the most popular. 
We will see later that this algorithm is a solver of \eqref{eq:mainproblem0}, with $\mathcal{M}=\mathcal{M}_{a,n}$.
Lately, explicit variational approaches \cite{schmaltz2010electrostatic,DeGoes} have been developed. The work of The work of de Goes \emph{et al} \cite{DeGoes} is closely related to our paper since they propose solving \eqref{eq:mainproblem0}, where $D$ is the $W_2$ transportation distance and $\mathcal{M}=\mathcal{M}_{f,n}$.
This sole problem is by no means limited to stippling and it is hard to provide a comprehensive list of applications. A few of them are listed in the introduction of \cite{xin2016centroidal}.

\paragraph{Best approximation with curves}

Another problem that is met frequently is to approximate a density by a curve. 
This can be used for non photorealitistic rendering of images or sculptures \cite{kaplan2005tsp,akleman2013hamiltonian}.
It can also be used to design trajectories of the nozzle of 3D printers \cite{chen2017line}.
It was also used for the generation of sampling schemes \cite{boyer2016generation}.

Apart from the last application, this problem is usually solved with methods that are not clearly expressed as an optimization problem.

\paragraph{Best approximation with arbitrary objects}

Problem \eqref{eq:mainproblem0} encompasses many other applications such as the optimization of networks \cite{gastner2006optimal}, texture rendering or non photorealistic rendering \cite{1210867,hiller2003beyond,schlechtweg2005renderbots,kim2009stippling,Jackbdu},  or sampling theory \cite{boyer2014algorithm}.

Overall, this paper unifies many problems that are often considered as distinct with specific methods.

\subsubsection{Numerical optimal transport}

In order to quantify the distance between the two measures, we use transportation distances \cite{monge1781memoire,kantorovich1942translocation,villani2003topics}.
In our work, we are interested mostly in the semi-discrete setting, where one measure is a density and the other is discrete. In this setting, the most intuitive way to introduce this distance is via Monge's transportation plan and allocation problems. Given an atomic measure 
$\nu\in \mathcal{M}_{\fdg{a,n}}$ and a  measure $\mu$ with density, a transport plan $T\in \mathcal{T}(\bx,\bw)$ is a mapping $T:\Omega\to \{\bx[1],\hdots, \bx[n]\}$ such that $\forall 1\leq i \leq n, \mu(T^{-1}(\bx[i]))=\bw[i]$. 
In words, the mass at any point $x\in \Omega$ is transported to point $T(x)$. 
In this setting, the $W_2$ transportation distance is defined by:
\begin{equation} \label{def:W2}
 W_2^2(\mu,\nu) = \inf_{T\in \mathcal{T}(\bx,\bw)} \int_{\Omega} \|x-T(x)\|_2^2  d \mu (x),
\end{equation}
and the minimizing mapping $T$ describes the optimal way to transfer $\mu$ to $\nu$.

Computing the transport plan $T$ and the distance $W_2$ is a challenging optimization problem.
In the semi-discrete setting, the paper \cite{aurenhammer1998minkowski} provided an efficient method based on ``power diagram'' or ``Laguerre diagram''. This framework was recently further improved and analyzed recently in  \cite{DeGoes,CGF:CGF2032,levy2015numerical,kitagawa2016newton}. 
The idea is to optimize a concave cost function with second-order algorithms. 
We will make use of those results in the paper, and improve them by stabilizing them while keeping the second-order information.

\subsubsection{Numerical projections on curve spaces}

Projecting curves on admissible sets is a basic brick for many algorithms. 
For instance, mobile robots are subject to kinematic constraints (speed and acceleration), while steel wire sculptures have geometric constraints (length, curvature). 

While the projection on kinematic constraints is quite easy, due to convexity of the underlying set \cite{chauffert2014gradient}, we believe that this is the first time projectors on sets defined through intrinsic geometry are designed. Similar ideas have been explored in the past. For instance, curve shortening with mean curvature motion \cite{evans1991motion} is a long-studied problem with multiple applications in computer graphics and image processing \cite{yezzi1998modified,moisan1998affine,tagliasacchi2012mean}.
The proposed algorithms allow exploring new problems such as curve lengthening with curvature constraints.

\subsection{Paper outline}

The rest of the paper is organized as follows. We first outline the overarching algorithm in Section~\ref{sec:numerical}. In Sections~\ref{sec:Wasserstein} and~\ref{sec:points}, we describe more precisely and study the theoretical guarantees of the algorithms used respectively for computing the Wasserstein distance, and for optimising the positions and weights of the points.  We describe the relationships with previous models in Section~\ref{sec:compare}. 
The algorithms in Sections~\ref{sec:Wasserstein} and~\ref{sec:points} are enough for, say, halftoning, but do not handle constraints on the points. In Section~\ref{sec:curves}, we add those constraints and design algorithms to make projections onto curves spaces with bounded speed and acceleration, or bounded length and curvature.
Finally some application examples and results are shown in Section~\ref{sec:applications}.


\section{The minimization framework}
\label{sec:numerical}

In this section, we show how to numerically solve the best approximation problem:
\begin{equation}\label{eq:continuousmeasureproblem}
 \inf_{\nu \in \mathcal{M}} W_2^2(\nu,\mu),
\end{equation}
where $\mathcal{M}$ is an arbitrary set of measures supported on $\Omega=[0,1]^2$.

\subsection{Discretization}

Problem \eqref{eq:continuousmeasureproblem} is infinite-dimensional and first needs to be discretized to be solved using a computer. We propose to approximate $\mathcal{M}$ by a subset $\mathcal{M}_n\subseteq\mathcal{M}_{a,n}$ of the atomic measures with $n$ atoms. The idea is to construct $\mathcal{M}_n$ as
\begin{equation}\label{eq:discretizedset}
 \mathcal{M}_n= \{ \nu(\bx,\bw), \bx\in \bX_n, \bw\in \bW_n\},
\end{equation}
where the mapping $\nu:(\Omega^n\times \Delta_{n-1})\to \mathcal{M}_{a,n}$ is defined by
\begin{align}
    \label{eq:nu}
    \nu(\bx,\bw) = \sum_{i=1}^n \bw[i] \delta_{\bx[i]}.
\end{align}
The constraint set $\bX_n\subseteq \Omega^n$  describes interactions between points and the set $\bW_n\subseteq \Delta_{n-1}$ describes the admissible weights.

 We have shown in \cite{chauffert2017projection} that for any subset $\mathcal{M}$ of the probability measures, it is possible to construct a sequence of approximation spaces $(\mathcal{M}_n)_{n\in \mathbb{N}}$ of the type \eqref{eq:discretizedset}, such that the solution sequence $(\nu_n^*)_{n\in \mathbb{N}}$ of the discretized problem
\begin{equation}\label{eq:discretizedmeasureprobleminitial}
 \inf_{\nu \in \mathcal{M}_n} W_2^2(\nu,\mu),
\end{equation}
converges weakly along a subsequence to a global minimizer $\nu^*$ of the original problem \eqref{eq:continuousmeasureproblem}. 
Let us give a simple example: assume that $\mathcal{M}$ is a set of pushforward measures of curves parameterized by a $1$-Lipchitz function on $[0,1]$. This curve can be discretized by a sum of $n$ Dirac masses with a distance between consecutive samples bounded by $1/n$. It can then be shown that this space $\mathcal{M}_n$ approximates $\mathcal{M}$ well, in the sense that each element of $\mathcal{M}_n$ can be approximated with a distance $O(1/n)$ by an element in $\mathcal{M}$ and vice-versa \cite{chauffert2017projection}. We will show explicit constructions of more complicated constraint sets $\bX_n$ and $\bW_n$ for measures supported on curves in Section~\ref{sec:curves}.

The discretized problem \eqref{eq:discretizedmeasureprobleminitial} can now be rewritten in a form convenient for numerical optimization:
\begin{equation}\label{eq:mainproblem}
 \min_{\bx \in \bX, \bw \in \bW} F(\bx,\bw),
\end{equation}
where we dropped the index $n$ to simplify the presentation and where 
\begin{equation}
 F(\bx,\bw) = \frac{1}{2} W_2^2(\nu(\bx,\bw),\mu).
\end{equation}

\subsection{Overall algorithm}

In order to solve \eqref{eq:mainproblem}, we propose to use an alternating minimization algorithm: the problem is minimized alternatively in $\bx$ with one iteration of a variable metric projected gradient descent and then in $\bw$ with a direct method. Algorithm~\ref{alg:mainalgo} describes the procedure in detail. 

A few remarks are in order. 
First notice that we are using a variable metric descent algorithm with a metric $\Sigma_k\succ 0$. 
Hence, we need to use a projector defined in this metric by:
\begin{align*}
    \Pi_\bX^{\Sigma_k}(\bx_0) & :=\Argmin_{\bx\in \bX} \|\bx-\bx_0\|^2_{\Sigma_k} \qquad \text{with} \\
    \|\bx-\bx_0\|^2_{\Sigma_k}          & \phantom{:}=\langle \Sigma_k (\bx-\bx_0),(\bx-\bx_0) \rangle.
\end{align*}
Second, notice that $\bX$ may be nonconvex. Hence, the projector $\Pi^{\Sigma_k}_\bX$ on $\bX$ might be a point-to-set mapping. In the $\bx$-step, the usual sign $=$ is therefore replaced by $\in$.

There are five major difficulties listed below to implement this algorithm:
\begin{description}
 \item[$\psi$ step: ]{How to compute efficiently $F(\bx,\bw)$?}
 \item[$\bw$ step: ]{How to compute $\displaystyle\argmin_{\bw \in \bW} F(\bx,\bw)$?}
 \item[$\bx$ step: ]{How to compute the gradients $\nabla_\bx F$ and the metric $\Sigma_k$?}
 \item[$\Pi$ step: ]How to implement the projector $\Pi^{\Sigma_k}_\bX$?
 \item[Generally: ]How to  accelerate the convergence given the specific problem structure?
\end{description}

The next four sections provide answers to these questions.

\begin{algorithm}
\caption{Alternating projected gradient descent to  minimize \eqref{eq:mainproblem0}. \label{alg:mainalgo}}
\begin{algorithmic}[1]
\Require  Oracle that computes $F$ \Comment{$\psi$-step}.
\Require Projectors $\Pi_\bX$ on $\bX$. 
\State \textbf{Inputs:} 
\State Initial guess $\bx_0$ 
\State Target measure $\mu$
\State Number of iterations $Nit$.
\State \textbf{Outputs:} 
\State An approximation $(\hat \bx,\hat \bw)$ of the solution of \eqref{eq:mainproblem0}.
\For{$k=0$ to $Nit-1$}
\State $\displaystyle  \bw_{k+1}=\argmin_{\bw \in \bW}(F(\bx_{k}, \bw)) $ \Comment{$\bw$-step}
\State Choose a positive definite matrix $\Sigma_k$, a step $s_k$. 
\State ${\bf y}_{k+1}=\bx_k - s_k \Sigma_k^{-1} \nabla_\bx F(\bx_k, \bw_{k+1})$. \Comment{$\bx$-step}
\State $\bx_{k+1}\in \Pi^{\Sigma_k}_\bX\left({\bf y}_{k+1}\right) \label{eq:PGDx}$ \Comment{$\Pi$-step}
\EndFor
\State Set $\hat \bx = \bx_{Nit}$ and $\hat \bw = \bw_{Nit}$.
\end{algorithmic}
\end{algorithm}
Note that if there are no constraints like in halftoning or stippling, there is no projection and the $\Pi$-step is trivial: $\bx_{k+1}={\bf y}_{k+1}$.

\section{Computing the Wassertein distance $F$ : $\psi$-step}
\label{sec:Wasserstein}
\subsection{Semi-discrete optimal transport}

In this paragraph, we review the main existing results about semi-discrete optimal transport and explain how it can be computed. 
Finally, we provide novel computation algorithms that are more efficient and robust than existing approaches.
We work under the following hypotheses.
\begin{assumption}\label{assumptions}\ 
\begin{itemize}
    \item The space $\Omega$ is a compact convex polyhedron, typically the hypercube.
 \item $\mu$ is an absolutely continuous probability density function w.r.t. the Lebesgue measure. 
 \item $\nu$ is an atomic probability measure supported on $n$ points. 
\end{itemize}
\end{assumption}

Let us begin by a theorem stating the uniqueness of the optimal transport plan, which is a special case of Theorem 10.41 in the book by \cite{villani2008optimal}.
\begin{theorem}
Under Assumption~\ref{assumptions}, there is a unique optimal transportation plan  $\mu$-a.e.$T^*$, solution of problem \eqref{def:W2}.
\end{theorem}

Before further describing the structure of the optimal transportation plan, let us introduce a fundamental tool from computational geometry  \cite{aurenhammer1991voronoi}.
\begin{definition}[Laguerre diagram]
Let $\bx \in \Omega^n$ denotes a set of point locations and $\bpsi\in \R^n$ denotes a weight vector.
 The Laguerre cell $\mathcal{L}_i$ is a closed convex polygon set defined as
\begin{equation}
    \label{eq:cells}
    \mathcal{L}_i(\bpsi,\bx) = \{ x \in \Omega, \forall 1\leq j \leq n, j\neq i,
    \|x-\bx[i]\|_2^2 - \bpsi[i] \leq \|x-\bx[j]\|_2^2 - \bpsi[j]  \}.
\end{equation}
\end{definition}
The Laguerre diagram generalizes the Voronoi diagram, since the latter is obtained by taking $\bpsi=0$ in equation \eqref{eq:cells}.

The set of Laguerre cells partitions $\Omega$ in polyhedral pieces.
It can be computed in $O(n\log(n))$ operations for points located in a plane \cite{aurenhammer1991voronoi}. 
In our numerical experiments, we make use of the  CGAL library to compute them \cite{cgal}. 
We are now ready to describe the structure of the optimal transportation plan $T^*$, see \cite[Example 1.9]{gangbo1996geometry}.
\begin{theorem}\label{thm:structureOT}
Under Assumption~\ref{assumptions}, there exists a vector $\bpsi^*\in \R^n$, such that 
\begin{equation}
  \fdg{(T^*)}^{-1}(\bx[i]) = \mathcal{L}_i(\bpsi^*,\bx).
\end{equation}
\end{theorem}
In words,  $\fdg{(T^*)}^{-1}(\bx[i])$ is the set where the mass located at point $\bx[i]$ is sent by the optimal transport plan. Theorem~\ref{thm:structureOT} states that this set is a convex polygon, namely the Laguerre cell of $\bx[i]$ in the tessellation with a weight vector $\bpsi^*$. 
More physical insight on the interpretation of $\bpsi^*$ can be found in \cite{levy2018notions}.
From a numerical point of view, the Theorem~\ref{thm:structureOT} allows transforming the infinite dimensional problem \eqref{def:W2} into the following finite dimensional problem:
\begin{equation}
    \label{eq:concave}
    W_2(\mu,\nu) =  \max_{\bpsi \in \mathbb{R}^n}  g(\bpsi,\bx,\bw),
\end{equation}
where 
\begin{equation}
\label{eq:defg}
g(\bpsi,\bx,\bw)  = \sum_{i=1}^n \int_{\mathcal{L}_i(\bpsi,\bx)} \left(\left\lVert \bx[i] - x  \right\rVert^2 - \bpsi[i]\right) \, \mathrm{d}\mu(x) + \sum_{i=1}^n \bpsi[i] \bw[i]. 
\end{equation}

Solving the problem \eqref{eq:concave} is the subject of numerous recent papers, and we suggest an original approach in the next section.

\subsection{Solving the dual problem}\label{sec:solvingdual}

In this paragraph, we focus on the resolution of \eqref{eq:concave}, i.e. computing the transportation distance numerically. 
The following proposition summarizes some concavity and differential properties of the functional $g$.  

\begin{proposition}\label{prop:derivatives}
Under Assumption~\ref{assumptions}, the function $g$ is concave with respect to the variable $\bpsi$ and differentiable with a Lipschitz gradient.
Its gradient is given by:
\begin{equation}
\frac{\partial g}{\partial \bpsi_i}=\bw[i]-\mu(\mathcal{L}_i(\bpsi,\bx)).
\end{equation}

In addition, if $\rho\in C^1(\Omega)$, then $g$ is also twice differentiable w.r.t. $\bpsi$ almost everywhere and - when defined - its second order derivatives are given by:
\begin{equation}\label{eq:HessianH}
 \frac{\partial^2 g}{\partial \bpsi_i\partial \bpsi_j}= \int_{\partial \mathcal{L}_i(\bpsi,\bx) \cap \partial \mathcal{L}_j(\bpsi,\bx)}\frac{d\mu(x)}{\Vert \bx[i]-\bx[j]\Vert} \text{ if } i\ne j. 
\end{equation}
The formula for the diagonal term $\frac{\partial^2 g}{\partial \bpsi_i\partial \bpsi_i}$ is given by the closure relationship
\begin{equation}
\forall 1\leq i \leq n, \quad \sum_{j=1}^n \frac{\partial^2 g}{\partial \bpsi_i\partial \bpsi_j} =0.
\end{equation}
\end{proposition}
\begin{proof}
 Most of these properties have been proved in \cite{kitagawa2016newton} and refined in \cite{Fred2018differentiation}. The Lipschitz continuity of the gradient seems to be novel.
 
\emph{Twice differentiability.} If a Laguerre cell $\mathcal{L}_i(\bpsi,\bx)$ is empty, it remains so for small variations of $\bpsi$, by the definition \eqref{eq:cells}.
It remains to prove that the set of $\bpsi$ for which there exists nonempty Laguerre cells with zero measure is negligible. The fact that $\mathcal{L}_i(\bpsi,\bx)$ is nonempty of zero Lebesgue measure means that it is either a segment or a point. We consider the case, where the points $\bx$ are in generic position, meaning that any three distinct points are not aligned. This implies that $\mathcal{L}_i(\bpsi,\bx)$ is a singleton $\{x\}$ since the boundaries of Laguerre cell cannot be parallel. We further assume that $x$ belongs to the interior of $\Omega$. Under those assumptions, $x$ necessarily satisfies at least $3$ equalities of the form 
\begin{equation}\label{eq:Li_singleton}
\|x-\bx[i]\|_2^2 - \bpsi[i] = \|x-\bx[j_k]\|_2^2 - \bpsi[j_k], 
\end{equation}
for some $j_k \neq i$ (i.e. $\mathcal{L}_i(\bpsi,\bx)$ is the intersection of at least $3$ half spaces). 
The set of $\bpsi$ allowing to satisfy a system of equations of the form \eqref{eq:Li_singleton} is of co-dimension at least $1$. Indeed, this system implies that $x$ is the intersection of 3 lines, each  perpendicular to one of the segment $[\bx[i],\bx[j_k]]$ and translated along the direction of this segment according to $\bpsi$. Now, by taking all the finitely many sets of quadruplets $(i,j_1,j_2,j_3)$, we conclude that the set of $\bpsi$ allowing to make $\mathcal{L}_i(\bpsi,\bx)$ a singleton is of zero Lebesgue measure. It remains to treat the case of $x$ belonging to the boundary of $\Omega$. This can be done similarly, by replacing one or more equalities in~\ref{eq:Li_singleton}, by the equations describing the boundary. 

The case of points in non generic position can also be treated similarly, since for $3$ aligned points at least $2$ equations of the form \eqref{eq:Li_singleton} allow to turn $\mathcal{L}_i(\bpsi,\bx)$ into a line of zero Lebesgue measure. A solution of such a system is also of co-dimension $1$. 

\emph{Lipschitz gradient.} In order to prove the Lipschitz continuity of the gradient, we first remark that the Laguerre cells defined in \eqref{eq:cells} are intersections of half spaces with a boundary that evolves linearly w.r.t. $\bpsi$. The rate of variation is bounded from above by $\eta = \frac{1}{\min_{i\neq j} \|\bx[i]-\bx[j]\|_2}$. Hence, letting $\Delta \in \R^n$ denote a unit vector, a rough bound on the variation of a single cell is:
\begin{equation*}
 \left\| \mu(\mathcal{L}_i(\bpsi+t\Delta,\bx)) - \mu(\mathcal{L}_i(\bpsi,\bx)) \right \| \leq t (n-1)\|\rho\|_\infty \eta \mathrm{diam}(\Omega).
\end{equation*}
Summing this inequality over all cells, we get that 
\begin{equation*}
 \| \nabla_{\bpsi} g(\bpsi+t\Delta,\bx) - \nabla_{\bpsi} g(\bpsi+t\Delta,\bx) \|_2 \leq t n^{3/2} \|\rho\|_\infty \eta \mathrm{diam}(\Omega).
\end{equation*}
Notice that this upper-bound is very pessimistic. For instance, applying Gersgorin's circle theorem shows that - when defined - the minimum eigenvalue of the Hessian matrix given in \eqref{eq:HessianH} is bounded below by $-n\eta \|\rho\|_\infty\mathrm{diam}(\Omega)$.
 
\end{proof}

In addition, the following proposition given in \cite[Thm. 6]{merigot2018algorithm} shows that the function $g$ is well behaved around the minimizers. 

\begin{proposition}\label{prop:regularity_minimizer}
If $\min_{x\in \Omega}\rho(x)>0$ and that the points $(\bx[i])$ are pairwise disjoint, Problem \eqref{eq:concave} admits a unique maximizer, up to the addition of constants. The function $g$ is twice differentiable in the vicinity of the minimizers and strongly concave on the set of vectors with zero mean.
\end{proposition}

Many methods have been proposed in the literature to compute the optimal vector $\bpsi^*$, with the latest references providing strong convergence guarantees \cite{aurenhammer1998minkowski,DeGoes,CGF:CGF2032,levy2015numerical,kitagawa2016newton}.
This may give the false impression that the problem has been fully resolved: in practice the conditions guaranteeing convergence are often unmet. For instance, it is well-known that the convergence of first-order methods depends strongly on the Lipschitz constant of the gradient \cite[Thm 2.1.7]{nesterov2013introductory}. Unfortunately, this Lipschitz constant may blow up depending on the geometry of the point set $\bx$ and the regularity of the density $\rho$, see Remark~\ref{rem:Lip}.
On the other hand, the second-order methods heavily depend on the H\"older regularity of $g$ \cite{jarre2016simple,grapiglia2017regularized}. Unfortunately, it can be shown that $g$ is H\"older with respect to $\psi$ only under certain circumstances. In particular, the mass of the Laguerre cells $\mu(\fdg{\mathcal L}_i(\bpsi,\bx))$ should not vanish \cite[Remark 4.2]{kitagawa2016newton}. Hence, only first-order methods should be used in the early steps of an optimization algorithm, and the initial guess should be well-chosen due to slow convergence. Then, second-order methods should be the method of choice.

The Levenberg-Marquardt algorithm and the trust-region methods \cite{wright1999numerical} are two popular solutions that interpolate between first- and second- order methods automatically. Unfortunately, to the best of our knowledge, the existing convergence theorems rely on a global $C^2$-regularity of the functional, which is not satisfied here. In this work, we therefore advocate the use of a regularized Newton method \cite{polyak2009regularized}, which retains the best of first and second order methods: a global convergence guarantee and a locally quadratic convergence rate. The algorithm reads as follows:
\begin{equation}\label{eq:regularized_Newton}
 \bpsi_{k+1} = \bpsi  - t_k (A(\bpsi) + \|\nabla_{\bpsi} g(\bpsi_k) \|_2 \mathrm{Id})^{-1}\nabla_{\bpsi} g(\bpsi_k),
\end{equation}
where 
\begin{equation}
 A(\bpsi)=\begin{cases}
           \nabla^2_{\bpsi}g(\bpsi_k) & \textrm{ if } \nabla^2_{\bpsi} g \textrm{ is defined at } \bpsi_k, \\
           0 &  \textrm{otherwise}.
          \end{cases}
\end{equation}
The algorithm is implemented on the set of vectors with zero mean to ensure the uniqueness of a solution, see Proposition~\ref{prop:regularity_minimizer}.

Without the term $\|\nabla_{\bpsi} g(\bpsi_k) \|_2 \mathrm{Id}$, the equation \eqref{eq:regularized_Newton} would simplify to a pure Newton algorithm. The addition of this term makes the method \eqref{eq:regularized_Newton} similar to a Levenberg-Marquardt algorithm, with the important difference that the regularization parameter is set automatically to $\|\nabla_{\bpsi} g(\bpsi_k) \|_2$. The rationale behind this choice is that the gradient vanishes close to the minimizer, making \eqref{eq:regularized_Newton} similar to a damped Newton method and that the gradient amplitude should be large far away from the minimizer, making \eqref{eq:regularized_Newton} closer to a pure gradient descent.

Following \cite{polyak2009regularized}, we get the following proposition.
\begin{proposition}
 Algorithm \eqref{eq:regularized_Newton} implemented with step-sizes $t_k$ satisfying the Wolfe conditions converges globally to the unique maximizer of \eqref{eq:concave}. In addition, the convergence is quadratic in the vicinity of the minimizer.
\end{proposition}
To the best of our knowledge, this is the first algorithm coming with a global convergence guarantee. Up to now, the convergence was only local \cite{kitagawa2016newton}. 

To further accelerate the convergence, the method can be initialized with the multi-scale approach suggested in \cite{CGF:CGF2032}. In practice, replacing \eqref{eq:regularized_Newton} by a standard Levenberg-Marquardt method i.e $\bpsi_{k+1} = \bpsi  - (A(\bpsi) + c_k \mathrm{Id})^{-1}\nabla_{\bpsi} g(\bpsi_k)$ yields similar rate of convergence. In this case $c_k > 0$ is interpreted as the ``step'' of the descent method and it is decreased or increased following Wolfe criterion.

\begin{remark}[High Lipschitz constant of the gradient]\label{rem:Lip}
  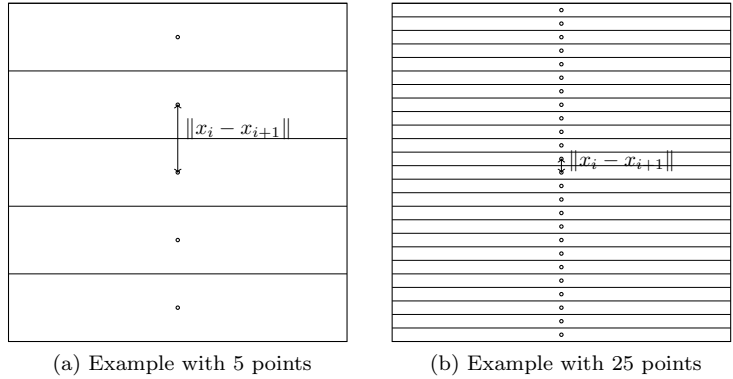
\begin{figure}
   \centering
    \subfloat[Example with 5 points]{\input{5pts}\label{fig:Lip5}} \hspace{1em}
    \subfloat[Example with 25 points]{\input{25pts}\label{fig:Lip25}}
    \caption{Configurations of points generating a high Lipschitz constant for the gradient of $g$ in $\bpsi$.\label{fig:Lip}}
  \end{figure}

  In this example illustrated by Figure~\ref{fig:Lip}, we show that the Lipschitz constant of the gradient can be arbitrarily large.
  We consider that $\mu$ is the uniform measure on $\Omega$ and that $\nu$ is an atomic measure supported on $n$ points aligned vertically and equispaced i.e. $ \bx [i] = \left(\frac{1}{2},\frac{1 + 2 i}{2n}\right) $ on $\Omega = [0,1]^2$.
 
  \fdg{In this case the Hessian is a multiple of the matrix of the $1d$ Laplacian with Neumann boundary conditions and the largest eigenvalue of $H$ scales as $4n$.}
 
  The Lipschitz constant hence blows up with the dimension since. 
  Notice that this situation is typical when it comes to approximating a density with a curve.
\end{remark}

\subsection{Numerical integration}\label{sec:numericalintegration}

The algorithm requires computing the integrals \eqref{eq:defg} and \eqref{eq:HessianH}. In all our numerical experiments, we use the following strategy. We first discretize the density $\rho$ associated to the target measure $\mu$ using a bilinear or a bi-cubic interpolation on a regular grid. Then, we observe that the volume integrals in Equation  \eqref{eq:defg} can be replaced by integrals of polynomials along the edges of the Laguerre diagram by using Green's formula. Hence computing the cost function, the Hessian or the gradient all boil down to computing edge integrals. 

Then, since the underlying density is piecewise polynomial, it is easy to see that only the \emph{first moments} of the measure $\mu$ along the edges are needed to compute all formula. We pre-evaluate the moments by using exact quadrature formulas and then use linear combinations of the moments to finish the evaluation.

To the best of our knowledge, this is a novel lightweight procedure. It significantly speeds up the calculations compared to former works \cite{CGF:CGF2032,DeGoes}, which enables discretization of the density $\rho$ over an arbitrary 3D mesh. After finishing this paper, we realized that the idea of using Green formulas was already suggested by \cite{xin2016centroidal}, although not implemented. It is to be noted that this idea is particularly well suited to Cartesian grid discretization of the target density $\rho$. Indeed, in this case, we take advantage of the fact that the intersection of the Laguerre cells and the grid can be computed analytically without search on the mesh. 

\section{Optimizing the weights and the positions: $\bw$ and $\bx$ steps}
\label{sec:points}
\subsection{Computing the optimal weights}

In this section, we focus on the numerical resolution of the following subproblem 
\begin{equation}\label{eq:optimW}
\argmin_{\bw \in \bW} F(\bx,\bw).  
\end{equation}

\subsubsection{Totally constrained $\bw$}
When $\bW=\{\bw\}$ is reduced to a singleton, the solution of \eqref{eq:optimW} is obviously given by $\bw^*=\bw$.

\subsubsection{Unconstrained minimization in $\bw$}

When $\bW$ is the simplex, the unconstrained minimization problem \eqref{eq:optimW} can be solved analytically.
\begin{proposition}
If $\ \bW=\Delta_{n-1}$, the solution $\bw^*$ of \eqref{eq:optimW} is given for all $1\leq i \leq n$ by
\begin{equation}
\bw^*[i]=\mu(\mathcal L_i(0,\bx)), 
\end{equation}
that is the volume (w.r.t. the measure $\mu$) of the $i$-th Laguerre cell with zero cost $\bpsi$, i.e. the $i$-th \emph{Vorono\"i cell}.
\end{proposition}
\begin{proof}
In expression \eqref{eq:defg}, the vector $\bpsi$ can be interpreted as a Lagrange multiplier for the constraint \[\mu(T^{-1}(\bx[i]))=\bw[i].\] Since the minimization in $\bw$ removes this constraint, the Lagrange multiplier might be set to zero. 

\end{proof}

\subsection{Gradient $\nabla_\bx F$ and the metric $\Sigma_k$}

The following proposition provides some regularity properties of $\nabla_\bx F$. It can be found in \cite{Fred2018differentiation}.
\begin{proposition}\label{prop:defgrad}
Let $\bpsi^*$  denote the maximizer of \eqref{eq:concave}. Assume that $\rho\in C^0\cap W^{1,1}(\Omega)$, that $\bw>0$, and that the points in $\bx$ are separated. Then $F$ is $C^2$ at $(\bx,\bw)$ with respect to the variable $\bx$ and its gradient $\nabla_\bx F(\bx,\bw)$ is given by the following formula:
\begin{equation}\label{eq:derivativex}
  \frac{\partial F(\bx,\bw)}{\partial \bx[i]} = \bw[i] \left( \bx[i] - \bb[i] \right) 
\end{equation}
where $\bb[i]$ is the barycenter of the $i$-th Laguerre cell $\mathcal{L}_i(\bpsi^*,\bx)$:
\begin{equation}\label{eq:Barycenter}
\bb[i]=\bb(\bx)[i] = \frac{\int_{\mathcal{L}_i(\bpsi^*,\bx)} x \mathrm{d}\mu(x)}{\int_{\mathcal{L}_i(\bpsi^*,\bx)}  \mathrm{d}\mu(x)}.
\end{equation}
\end{proposition}

Now, we discuss the choice of the metric $(\Sigma_k)$ in Algorithm~\ref{alg:mainalgo}. In what follows, we refer to the ``unconstrained case'' as the case where there is no $\Pi$-step in Algorithm~\ref{alg:mainalgo}. The metric used in our paper is the following:
\begin{equation}\label{eq:MetricSk}
 \Sigma_k=\diag( \mu(\mathcal{L}_i(\bpsi^*_k,\bx^*_k)))_{1\leq i \leq n}).
\end{equation}

We detail the rationale behind this choice below.
First, with the choice \eqref{eq:MetricSk}, we have $\bx_k  - \Sigma_k^{-1}\nabla_\bx F(\bx_k) = \bb(\bx_k)$. In the unconstrained case, this particular choice of $\Sigma_k$ amounts to moving the points $\bx$ towards their barycenters, which is the celebrated Lloyd's algorithm.

Beside this nice geometrical intuition, in the unconstrained case, the choice \eqref{eq:MetricSk} leads to an alternating direction minimization algorithm. Indeed, given a set of points, this algorithm computes the optimal transport plan ($\bpsi$-step). Then, fixing this transport plan and the associated Laguerre tessellation, the mass of the point is moved to the barycenter of the Laguerre cell, which is the optimal position for a given tessellation. This algorithm is widespread because it does not require additional line-search.

Third, in the unconstrained case, the choice \eqref{eq:MetricSk} leads to an interesting regularity property around the critical points. Assume that $\nabla_\bx F(\bx^*)=0$, i.e. that $\bx^*[i]=\bb(\bx^*)[i]$ for all $i$, then the mapping $\bx\mapsto\bb(\bx)$ is locally 1-Lipschitz \cite[Prop. 6.3]{du1999centroidal}. This property suggests that a variable metric gradient descent with metric $\Sigma_k$ and step size $1$ may perform well in practice for $\bX=\Omega^n$, at least around critical points.

Fourth, this metric is the diagonal matrix with coefficients obtained by summing the coefficients of the corresponding line of $H_{\bx \bx}[F]$, the Hessian of $F$ with respect to $\bx$ see \cite{Fred2018differentiation}. In this sense, $\Sigma_k$ is an approximation of $H_{\bx \bx}[F]$. In the unconstrained case Algorithm \ref{alg:mainalgo} can be interpreted as a quasi-Newton algorithm. 

A safe choice of the step $s_k$ in Algorithm \ref{alg:mainalgo} to ensure convergence could be driven by Wolfe conditions. In view of all the above remarks, it is tempting to use a gradient descent with the choice $s_k=1$. In practice, it gives a satisfactory rate of convergence for Algorithm \ref{alg:mainalgo}. For all the experiments presented in this paper, we therefore make this empirical choice. We provide some elements to justify the local convergence under a more conservative choice of parameters in Section \ref{sec:app_convergence}.

\section{Links with other models}
\label{sec:compare}
\subsection{Special cases of the framework}
\subsubsection{Lloyd's algorithm}

Lloyd's algorithm \cite{lloyd1982least} is well-known to be a specific solver for problem \eqref{eq:continuousmeasureproblem}, with $\bX=\Omega$ and $\bW=\Delta_{n-1}$, i.e. to solve the quantization problem with variable weights. We refer to the excellent review by \cite{du1999centroidal} for more details. It is easy to check that Lloyd's algorithm is just a special case of Algorithm~\ref{alg:mainalgo}, with the specific choice of metric
\begin{equation}
 \Sigma_k = \diag\left( \mu(\mathcal{L}_i(0,\bx)) \right).
\end{equation}

\subsubsection{Blue noise through optimal transport}
In, \cite{DeGoes}, the authors has proposed to perform stippling by using optimal transport distance. 
This application can be cast as a special case of problem \eqref{eq:continuousmeasureproblem}, with $\bX=\Omega$ and $\bW=\left\{\frac{\one}{n}\right\}$. 
The algorithm proposed therein is also a special case of algorithm~\ref{alg:mainalgo} with 
\begin{equation}
 \Sigma_k = \diag\left( \mu(\mathcal{L}_i(\phi^\star(\bx),\bx)) \right)=\frac{1}{n}
\end{equation}
and the step-size $\tau_k$ is optimized through a line search. Note however the extra cost of applying a line-search might not worth the effort, since a single function evaluation requires solving the dual problem \eqref{eq:concave}.

 \subsection{Comparison with electrostatic halftoning}
 
In \cite{schmaltz2010electrostatic,teuber2011dithering,fornasier2013consistency,chauffert2017projection}, an alternative to the $W_2$ distance was proposed, implemented and studied. Namely, the distance $D$ in \eqref{eq:mainproblem0} is defined by
\begin{equation}\label{eq:distanceconvolution}
 D(\nu,\mu)=\frac{1}{2}\|h\star (\nu-\mu)\|^2_{L^2(\Omega)},
\end{equation}
where $h$ is a smooth convolution kernel and $\star$ denotes the convolution product. 
This distance can be interpreted intuitively as follows: the measures are first blurred by a regularizing kernel to map them in $L^2(\Omega)$ and then compared using a simple $L^2$ distance.
It appears in the literature under different names such as Maximum Mean Discrepancies, kernel norms or blurred SSD.
In some cases, the two approaches are actually quite similar from a theoretical point of view. Indeed, it can be shown that the two distances are strongly equivalent under certain assumptions \cite{peyre2016comparison}.

The two approaches however differ significantly from a numerical point of view. Table~\ref{tab:comparison} provides a quick summary of the differences between the two approaches. We detail this table below.
\begin{itemize}
 \item The theory of optimization is significantly harder in the case of optimal transport since it is based on a subtle mix between first and second order methods. 
 \item The convolution-based algorithms require the use of methods from applied harmonic analysis dedicated to particle simulations such as fast multiple methods (FMM) \cite{greengard1987fast} or non-uniform Fast Fourier Transforms (NUFFT) \cite{potts2003fast}. On their side, the optimal transport based approaches require the use of computational geometry tools such as Voronoi or Laguerre diagrams. The former has been parallelized efficiently on CPU and GPU and turnkey toolboxes are now available, while the latter seem to be less accessible for now and some implementations are intrinsically serial.
 \item The examples provided here are only two-dimensional. Many applications in computer graphics require dealing with 3D problems or larger dimensional problems (e.g. clustering problems). In that case, the numerical complexity of convolution based problems seems much better controlled: it is only linear in the dimension $d$ (i.e. $O(dn \log(n))$), while the exact computation of Laguerre diagrams requires  in average $O(n^{\lceil d/2 \rceil} )$ operations, the worst case time complexity for $d=2$ is $O(n\log n)$\cite{aurenhammer1987power}. Hence, for a large number of particles, the approach suggested here is mostly restricted to $d=2$.
 \item In terms of computational speed for 2D applications, we observed that the optimal transport based approach was usually between 1 and 2 orders of magnitude faster. This is mostly due to the fact that the descent algorithm based on optimal transport converges in significantly less iterations than that based on convolution distances.
 \item Finally, we did not observe significant differences in terms of approximation quality from a perceptual point of view.
\end{itemize}

  \begin{table}
  \centering
  \begin{tabular}{|l|l|l|}
  \hline
   & \textbf{Convolution} & \textbf{Optimal transport} \\
  \hline
  \textbf{Optimization} & \green{1st order}   & \red{Mix of 1st and 2nd} \\
  \hline
  \textbf{Computation} & FMM/NUFFT   & Power diagram \\
  \hline
  \textbf{Scaling to $d$} & \green{Linear}    & \red{Exponential} \\
  \hline
  \textbf{Speed in 2d} & \red{Slower} & \green{Faster} \\
  \hline 
  \textbf{Quality} & \green{Good} & \green{Good} \\
  \hline
  \end{tabular}
  \caption{A comparison between convolution and optimal transport based approximation of measures. \label{tab:comparison}}
  \end{table}

\subsubsection{Benchmark with other methods}
In this section we provide compare 4 methods: two versions of electrostatic halftoning \cite{chauffert2017projection}, ibnot (a semi-discrete optimal transport toolbox \cite{DeGoes}) and the code presented in this paper.

\paragraph{Choice of a stopping criterion}
  The comparison of different methods yields the question of a stopping criterion. Following \cite{schmaltz2010electrostatic}, the signal-to-noise ratio (SNR) of the original image and the stippled image convolved with a Gaussian function has been chosen. The standard deviation of the Gaussian is chosen as $1/\sqrt{n}$, which is the typical distance between points. Figure~\ref{fig:stopCrit} shows different values of this criterion for an increasing quality of stippling. In all the forthcoming benchmarks, the algorithms have been stopped when the underlying measured reached a quality of $31$dB.
  
  \def\scal{.24}
  \begin{figure}[!h]
   \centering
   \includegraphics[width=\scal \textwidth]{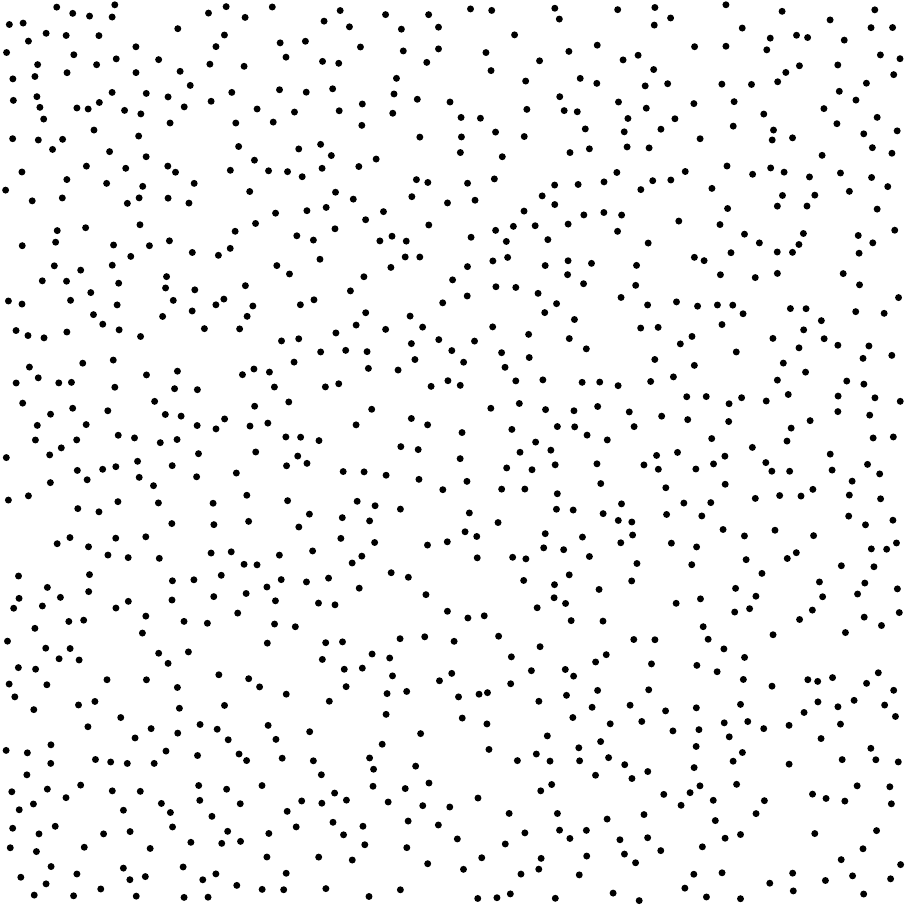}
   \includegraphics[width=\scal \textwidth]{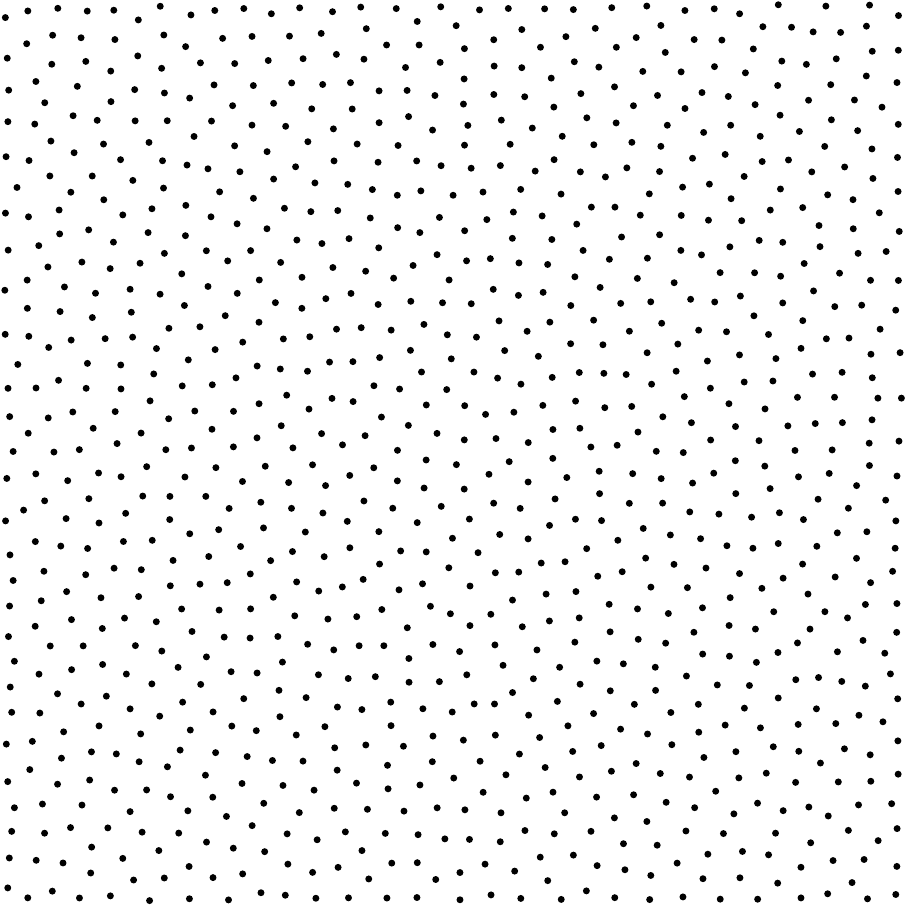}
   \includegraphics[width=\scal \textwidth]{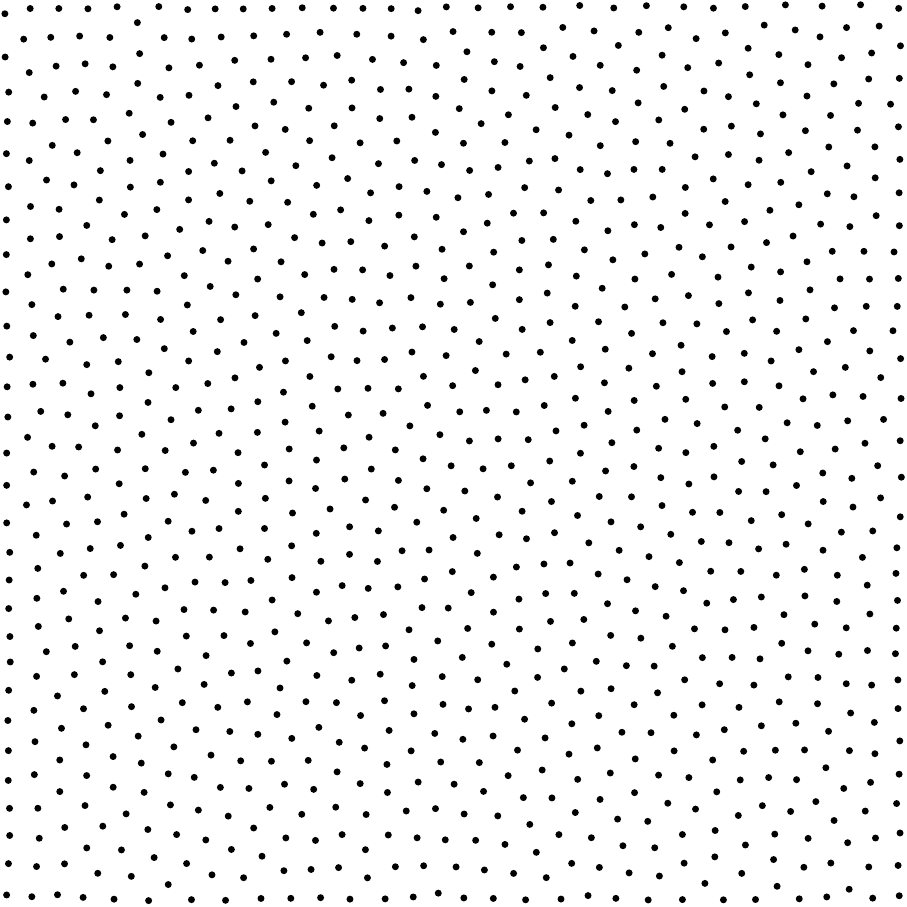}
   \includegraphics[width=\scal \textwidth]{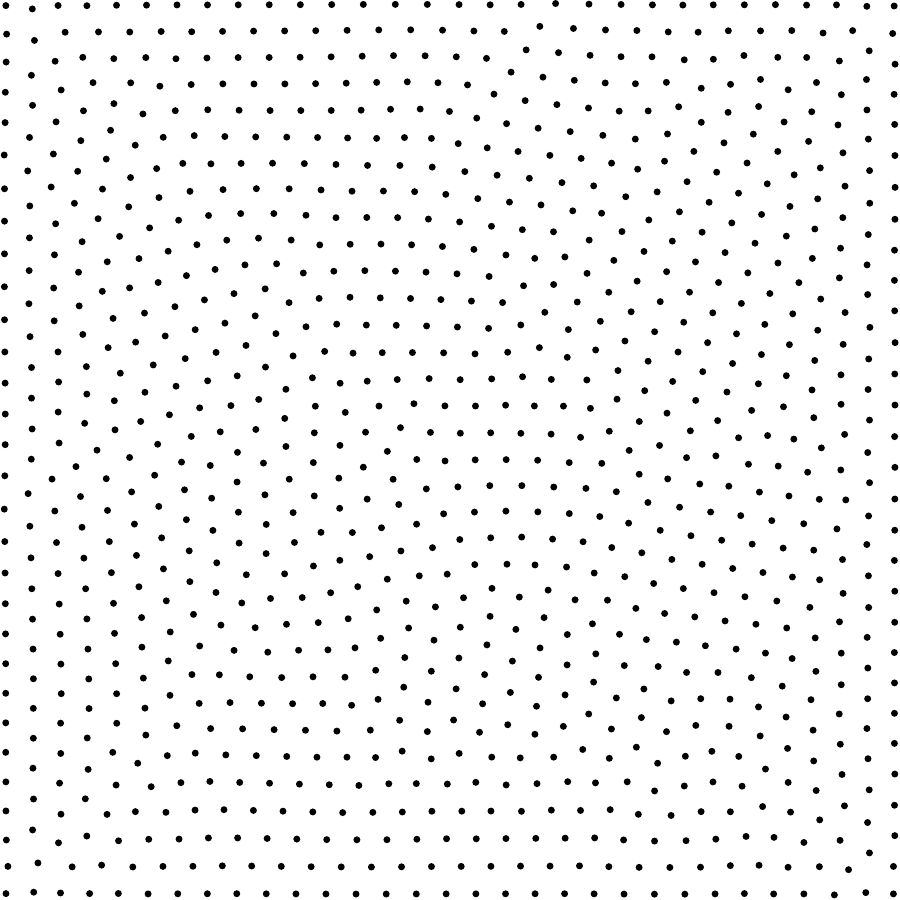}
   \caption{Evolution of the PSNR through iteration : top left $8$dB, top right $25$dB, bottom left $31$dB, bottom right $34$dB}\label{fig:stopCrit}

  \end{figure}

 \paragraph{Benchmarks}
The first benchmark  is illustrated in Table~\ref{tab:benchmark}. For this test the background measure has a constant density and consists of 1024$\times$1024 pixels. The number of Dirac masses increases from $2^{10}$ to $2^{18}$. The initialization is obtained by a uniform Poisson point process.
 
 \begin{table}[!ht]
  \centering
  \begin{tabular}{|l|c|c|c|c|}
   \hline
    \textbf{\# pts} & \textbf{Electro 20 cores} & \textbf{Electro BB 20 cores} & \textbf{ibnot 1 core} & \textbf{our 1 core}\\
   \hline
   \small{$2^{10}$} & $130.3$ | 317 & $34.4$ | 84 & $131.47$ | 15 & $4.03$ | 19\\
   \hline
   \small{$2^{11}$} & $240.1$ | 558 & $49.5$ | 115 & $165.42$ | 18 & $4.83$ | 23\\
   \hline
   \small{$2^{12}$} & $293.9$ | 637 & $47.8$ | 104 & $267.59$ | 22 & $10.86$| 19 \\
   \hline
   \small{$2^{13}$} & $415.2$ | 798 & $78.8$ | 152 & $235.87$ | 15 & $26.43$ | 20\\
   \hline
   \small{$2^{14}$} & $783.5$ | 1306 & $106.3$ | 177 & $344.77$ | 17 & $47.90$| 21\\
   \hline
   \small{$2^{15}$} & $1160.5$ | 1319 & $156.8$ | 178 & $598.60$ | 18 & $99.63$| 18\\
   \hline
   \small{$2^{16}$} & $4568.5$ | 3286  & $569.2$ | 410 & $1208.45$ | 20 & $252.24$| 26 \\
   \hline
   \small{$2^{17}$} & $15875.0$ | 4628 & $1676.3$ | 489 & $2498.58$ | 19& $620.90$| 19 \\
   \hline
   \small{$2^{18}$} & $TL$ &  $12125.2$ | 1103 & $5633.68$ | 23& $1136.51$| 21\\
   \hline
  \end{tabular}
  
  \caption{Times in second and number of iterations to achieve convergence for the uniform background measure. TL means that the computing time was too long and that we stopped the experiment before reaching the stopping criterion. For the Electrostatic halftoning, two different methods are compared : a gradient descent with constant step-size (2nd column), and the Brazalai Borwein method (third column).}\label{tab:benchmark}
 \end{table}
 
 In this case ibnot and our code present the same complexity and roughly the same number of $\bx$-step to achieve convergence. The time per iterations is significantly smaller in our code due to the use of Green's formula for integration (see Section~\ref{sec:numericalintegration}), which reduces the integration's complexity from $n$ to $\sqrt{n}$ where $n$ the number of pixels.
 We tested multiple versions of electrostatic halftoning, differing in the choice of the optimization algorithm. The code Electro 20 cores corresponds to a constant step-size gradient descent as proposed in \cite{schmaltz2010electrostatic}. The code Electro BB 20 cores is a gradient descent with a Barzilai-Borwein step-size rule. In our experience, this turns out to be the most efficient solver (e.g. more efficient than a L-BFGS algorithm). The algorithms have been parallelized with Open-MP and evaluated on a 20 cores machine using the NFFT to evaluate the sums \cite{potts2003fast}. 
 
 In Table~\ref{tab:benchmark}, the electrostatic halftoning algorithms are always slower than our optimal transport algorithms despite being multi-threaded on $20$ cores.
 
 The second test is displayed in Table~\ref{tab:benchmark2}. It consists in trying to approximate a non-constant density $\rho(x,y)$ equal to $2$ if $x<0.5$ and $0$ if $x>0.5$. In this test we also start from a uniform point process.
 
 \begin{table}[!ht]
  \centering
  \begin{tabular}{|l|c|c|c|}
   \hline
   \textbf{\# pts} & \textbf{Electro  BB 20 cores} & \textbf{ibnot 1 core} & \textbf{our 1 core}\\
   \hline
   \small{$2^{10}$} & $40.2$ | 99 & NC & $4.34$ | 24\\
   \hline
   \small{$2^{11}$} & $59.6$ | 142 & NC & $9.70$ | 21 \\
   \hline
   \small{$2^{12}$} & $84.5$ | 190 & NC & $15.36$| 19\\
   \hline
   \small{$2^{13}$} & $125.1$ | 248 & NC & $29.76$| 25 \\
   \hline
   \small{$2^{14}$} & $177.7$ | 282 & NC & $79.73$| 27 \\
   \hline
   \small{$2^{15}$} & $287.6$ | 298 & NC & $113.56$ |18\\
   \hline
   \small{$2^{16}$} & $746.3$ | 424 & NC & $280.02$| 17\\
   \hline
   \small{$2^{17}$} & $3052.9$ | 712 & NC & $703.39$| 21\\
   \hline
   \small{$2^{18}$} & $39546.1$ | 2022 & NC & $1315.01$|24\\
   \hline
  \end{tabular}
  \caption{Times in seconds and number of iterations to achieve convergence for the non uniform setting. NC stands for does not converge. TL stands for too long (exceeding 4 hours of computations on 20 cores).}\label{tab:benchmark2}
 \end{table}
 
 In this test ibnot always fails to converge since the Hessian in $\bpsi$ is not definite. In comparison with the previous test, our code is slightly slower since the first $\bpsi$-step requires to find a $\bpsi$ that maps Laguerre cells far away from the localization of their site ($\bx$). After the first $\bx$-step, the position of the site is decent and the optimization routine performs as well as the previous example. Again, the computing times of electro-static halftoning is significantly worse than the one of optimal transport. Notice in particular how the number of iterations needed to reach the stopping criterion increases with the number of points, while it remains about constant for the optimal transport algorithm.


\section{Projections on curves spaces}\label{sec:curves}

In this section, we detail a numerical algorithm to evaluate the projector $\Pi_{\bX}$, for spaces of curves with kinematic or geometric constraints.

\subsection{Discrete curves}

A discrete curve is a set of points $\bx\in \Omega^n$ with constraints on the distance between successive points. 
Let 
\begin{align*}
A_1^a&:\bx \to \begin{pmatrix}
               \bx[2]-\bx[1] \\ \vdots \\ \bx[n]-\bx[n-1] \\ \bx[1]-\bx[n]
              \end{pmatrix} \\
              \end{align*}
              and
              \begin{align*}
A_1^b&:\bx \to \begin{pmatrix}
               \bx[2]-\bx[1] \\ \bx[3]-\bx[2] \\ \vdots \\ \bx[n]-\bx[n-1]
              \end{pmatrix}
\end{align*}
denote the discrete first order derivatives operators with or without circular boundary conditions.
From hereon, we let $A_1$ denote any of the two operators. 
In order to control the distance between two neighboring points, we will consider two types of constraints: kinematic ones and geometrical ones.

\subsubsection{Kinematic constraints}

Kinematic constraints typically apply to vehicles: a car for instance has a bounded speed and acceleration. 
Bounded speeds can be encoded through inequalities of type
\begin{equation}\label{eq:boundedspeed}
 \|(A_1\bx)[i]\|_2 \leq \alpha_1, \forall i.
\end{equation}
Similarly, by letting $A_2$ denote a discrete second order derivative, which can for instance be defined by $A_2=A_1^TA_1$, we may enforce bounded acceleration through
\begin{equation}
 \|(A_2\bx)[i]\|_2 \leq \alpha_2, \forall i.
\end{equation}
The set $\bX$ is then defined by
\begin{equation}\label{eq:kinematicset}
 \bX=\{\bx\in \Omega^n, \|A_1\bx\|_{\infty,2} \leq \alpha_1, \|A_2\bx\|_{\infty,2} \leq \alpha_2\},
\end{equation}
where, for $\by=(\by[1],\hdots,\by[n])$, $\|\by\|_{\infty,p}=\sup_{1\leq i \leq n} \|\by[i]\|_p$.

\subsubsection{Geometrical constraints}

Geometrical constraints refer to intrinsic features of a curve such as its length or curvature. 
In order to control those quantities using differential operators, we need to parameterize the curve with its arc length.
Let $s:[0,T]\to \R^2$ denote a $C^2$ curve with arc length parameterization, i.e. $\|\dot s(t)\|_2=1, \forall t\in [0,T]$. 
Its length is then equal to $T$. Its curvature at time $t\in [0,T]$ is equal to $\kappa(t)=\|\ddot s(t)\|_2$.

In the discrete setting, constant speed parameterization can be enforced by imposing
\begin{equation}
\label{eq:constantspeed}
 \|(A_1\bx)[i]\|_2 = \alpha_1, \forall i.
\end{equation}
The total length of the discrete curve is then equal to $(n-1)\alpha_1$.

Similarly, when \eqref{eq:constantspeed} is satisfied, discrete curvature constraints can be captured by inequalities of type
\begin{equation}\label{eq:curvature}
 \|(A_2\bx)[i]\|_2 \leq \alpha_2, \forall i.
\end{equation}
Indeed, at a index $2\leq i \leq n-1$, we get:
\begin{align*}
 \|(A_2\bx)[i]\|_2^2 & = \|(\bx[i]- \bx[i-1]) - (\bx[i+1] - \bx[i]) \|_2^2 \\
 &=\|\bx[i]- \bx[i-1]\|_2^2 + \|\bx[i+1] - \bx[i] \|_2^2 \\
 &- 2\langle \bx[i]- \bx[i-1], \bx[i+1] - \bx[i]\rangle \\
 &= 2 \alpha_1^2 (1- \cos\left( \theta_i\right)),
\end{align*}
where $\theta_i=\angle \left(\bx[i]- \bx[i-1], \bx[i+1] - \bx[i] \right)$ is the angle between successive segments of the curve.
Hence, by imposing \eqref{eq:constantspeed} and \eqref{eq:curvature}, the angle $\theta_i$ satisfies
\begin{equation}
 |\theta_i| \leq \arccos\left( 1 - \frac{\alpha_2^2}{2\alpha_1^2} \right).
\end{equation}

In order to fix the length and bound the curvature, we may thus choose the set $\bX$ as
\begin{equation}\label{eq:geometricset}
 \bX=\{\bx\in \Omega^n, \|(A_1\bx)[i]\|_{2} = \alpha_1, \|A_2\bx\|_{\infty,2} \leq \alpha_2\}.
\end{equation}
Let us note already that this set is nonconvex, while \eqref{eq:kinematicset} was convex.

\subsubsection{Additional linear constraints}

In applications, it may be necessary to impose other constraints such as passing at a specific location at a given time, closing the curve with $x_1=x_n$ or having a specified mean value.
All those constraints are of form
\begin{equation}
 B\bx = \bb,
\end{equation}
where $B\in \R^{p\times 2n}$ and $\bb\in \R^p$ are a matrix and vector describing the $p$ linear constraints.

\subsubsection{Summary}

In this paper, we will consider discrete spaces of curves $\bX$ defined as follows:
\begin{equation}
 \bX = \{\bx \mbox{ such that } A_i\bx\in \bY_i, 1\leq i\leq m, B\bx =b\},
\end{equation}
The operators $A_i$ may be arbitrary, but in this paper, we will focus on differential operators of different orders.
The set $\bY_i$ describes the admissible set for the $i$-th constraint. 
For instance, to impose a bounded speed \eqref{eq:boundedspeed}, we may choose
\begin{equation}
 \bY_1=\{\by \in \R^{n\times 2}, \|\by_i\|_2\leq \alpha_1, \forall i\}.
\end{equation}
In all the paper, the set of admissible weights $\bW$ will be either the constant $\{\one/n\}$ or the canonical simplex $\Delta_{n-1}$.

\subsection{Numerical projectors}

The Euclidean projector $\Pi_\bX:\R^n\to\bX$ is defined for all $\bz\in \Omega^n$ by
\begin{align}
 \Pi_\bX(\bz) &= \Argmin_{\bx \in \bX} \frac{1}{2} \|\bx-\bz\|_2^2 \nonumber \\
 &= \Argmin_{\substack{A_k \bx \in \bY_k, 1\leq k\leq m \\ B\bx=\bb}} \frac{1}{2} \|\bx-\bz\|_2^2 \label{eq:projector}
\end{align}
When $\bX$ is convex, $\Pi_\bX(\bz)$ is a singleton. 
When it is not, there exists $\bz$ such that $\Pi_\bX(\bz)$ contains more than one element.
The objective of this section is to design an algorithm to find critical points of \eqref{eq:projector}. 

The specific structure of \eqref{eq:projector} suggests using splitting based methods \cite{combettes2011proximal}, able to deal with multiple constraints and linear operators. The sparse structure of differential operator makes the Alternating Direction Method of Multipliers (ADMM, \cite{glowinski2014alternating}), particularly suited for this problem. Let us turn \eqref{eq:projector} into a form suitable for the ADMM. 

Let $\gamma_1,\hdots,\gamma_m$ denote positive reals used as preconditionners. Define 
\begin{equation}
 A=\begin{pmatrix}
    \gamma _1 A_1 \\ \vdots \\ \gamma_m A_m 
   \end{pmatrix}
   \mbox{,} \quad
   \by=\begin{pmatrix}
	\by_1 \\ \vdots \\ \by_m
       \end{pmatrix}
\end{equation}
and 
\begin{equation}
 \bY= \gamma_1\bY_1\times \hdots\times \gamma_m \bY_m.
\end{equation}
Problem \eqref{eq:projector} then becomes
\begin{align}
 \Pi_\bX(\bz) &= \Argmin_{\substack{B\bx=\bb \\ A \bx = \by \\ \by \in \bY}} \frac{1}{2} \|\bx-\bz\|_2^2 \nonumber\\
 &= \Argmin_{A \bx = \by} f_1(\bx) + f_2(\by),\label{eq:ADMMsuited}
\end{align}
where $f_1(\bx) = \frac{1}{2}\|\bx-\bz\|_2^2+\iota_\bL(\bx)$, $f_2(\by)=\iota_\bY(\by)$, $\bL=\{\bx, B\bx=\bb\}$ denotes the set of linear constraints and the indicator $\iota_\bY$ of $\bY$ is defined by:
\begin{equation}
 \iota_\bY(\by)=\begin{cases}
           0 & \textrm{if } \by\in\bY, \\
           +\infty  & \textrm{otherwise}.
          \end{cases}
\end{equation}
The ADMM for solving \eqref{eq:ADMMsuited} is given in Algorithm \ref{alg:ADMMgeneric}.
Specialized to our problem, it yields Algorithm \ref{alg:ADMM}. The linear system can be solved with a linear conjugate gradient descent.

\begin{algorithm}
\caption{Generic ADMM. \label{alg:ADMMgeneric}}
\begin{algorithmic}[1]
\Statex \textbf{Inputs:} 
\Statex functions $f_1$ and $f_2$, matrix $A$, initial guess $(\bx_0,\blambda_0)$, parameter $\beta>0$.
\While{Stopping criterion not met}
\Statex $\displaystyle \by_{k+1}=\Argmin_{\by} f_2(\by) + \frac{\beta}{2}\|A\bx -\by_{k+1} + \blambda_k\|_2^2$.
\Statex $\displaystyle \bx_{k+1}=\Argmin_{\bx} f_1(\bx) + \frac{\beta}{2}\|A\bx -\by_{k+1} + \blambda_k\|_2^2$. 
\Statex $\blambda_{k+1}= \blambda_{k} + A\bx_{k+1} -\by_{k+1}$.
\EndWhile
\end{algorithmic}
\end{algorithm}

\begin{algorithm}
\caption{ADMM to solve the projection problem. \label{alg:ADMM}}
\begin{algorithmic}[1]
\Statex \textbf{Inputs:} 
\Statex Vector to project $\bz$, initial guess $(\bx_0,\blambda_0)$, matrices $A$ and $B$, projector $(\Pi_{\bY})$, $\beta>0$.
\While{Stopping criterion not met}
\Statex $\displaystyle \by_{k+1}= \Pi_{\bY}(A\bx_{k} + \blambda_k)$.
\Statex Solve \begin{align*} \begin{bmatrix}
                              \beta A^TA + I & B^T \\ B & 0 
\end{bmatrix} \begin{pmatrix} \bx_{k+1} \\ \bmu\end{pmatrix} 
&= \begin{pmatrix} \beta A^T(\by_{k+1}- \blambda_k)+ \bz \\ b\end{pmatrix}.\end{align*}
\Statex $\blambda_{k+1}= \blambda_{k} + A\bx_{k+1} -\by_{k+1}.$
\EndWhile
\end{algorithmic}
\end{algorithm}

\paragraph{Convergence issues}

The convergence and rate of convergence of the ADMM is a complex issue that depends on the properties of functions $f_1$ and $f_2$ and on the linear transform $A$. In the convex setting~\eqref{eq:kinematicset}, the sequence $(\bx_k)_k$ converges to $\Pi_{\bX}(\bz)$ linearly (see Corollary 2 in \cite{giselsson2017linear}). The behavior in a nonconvex setting~\eqref{eq:geometricset} is still mostly open despite recent advances  in \cite{li2015global}.
Nevertheless, we report that we observed convergence empirically towards critical points of Problem \eqref{eq:projector}. 

\paragraph{Choosing the coefficients $\beta$ and  $(\gamma_i)$}

Despite recent advances \cite{nishihara2015general}, a theory to select good values of $\beta$ and  $(\gamma_i)$ still seems lacking. 
In this paper, we simply set $\gamma_i=\|A_i\|_2$, the spectral norm of $A_i$.
In practice, it turns out that this choice leads to stable results. 
The parameter $\beta$ is set manually to obtain a good empirical behavior.
Notice that for a given application, it can be tuned once for all.

\subsection{Numerical examples}

\begin{figure*}[ht]
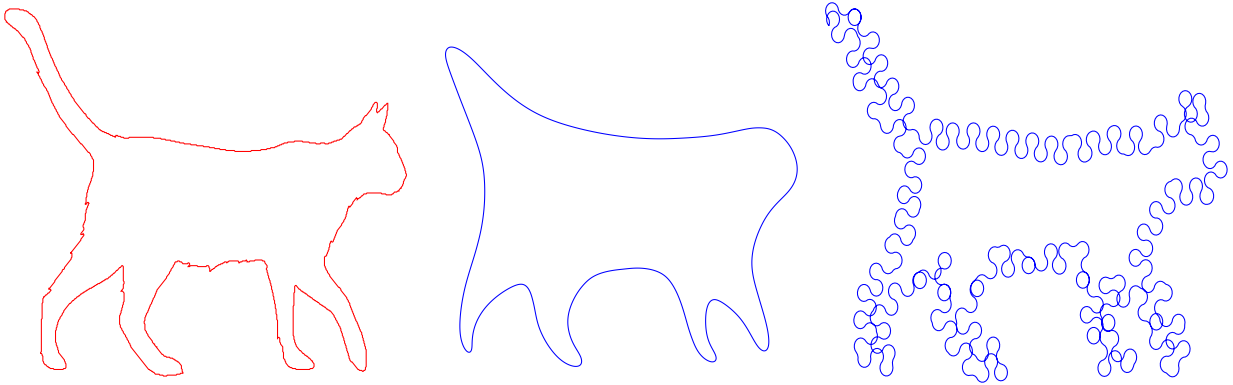

\centering
 \input{Cat_ini.tex}
 \input{Cat_8p00e-01_2p00e+00.tex}
 \input{Cat_3p00e+00_1p00e+00.tex}
 \caption{Examples of projections of a curve (in red) on spaces of curves with constraints (in blue). Center: projection on sets of curves with smaller length and bounded curvature. Right: projection on sets of curves with longer length and bounded curvature.\label{fig:CurveProjection}}
\end{figure*}

To illustrate the proposed method, we project the silhouette of a cat onto spaces of curves with fixed length and bounded curvature in Fig.~\ref{fig:CurveProjection}. In the middle, we see how the algorithm simplifies the curve by making it smaller and smoother. On the right, we see how the method is able to make the curve longer, by adding loops of bounded curvature, while still keeping the initial cat's shape.

\subsection{Multi-resolution implementation}

When $\bX$ is a set of curves, the solution of \eqref{eq:mainproblem} can be found more efficiently by using a multi-resolution approach. 
Instead of optimizing all the points simultaneously, it is possible to only optimize a down-sampled curve, allowing to get cheap warm start initialization for the next resolution. 

In our implementation, we use a dyadic scaling. We up-sample the curve by adding mid-points in between consecutive samples. The weights from one resolution to the next are simply divided by a factor of $2$.


\section{Applications}
\label{sec:applications}

\subsection{Non Photorealistic Rendering with curves}

In the following subsections we exhibit a few rendering results of images using different types of measures sets $\mathcal{M}$.

\subsubsection{Gray-scale images}

\begin{figure*}
  \centering
  \subfloat[Original]{{\includegraphics[width=.23\linewidth]{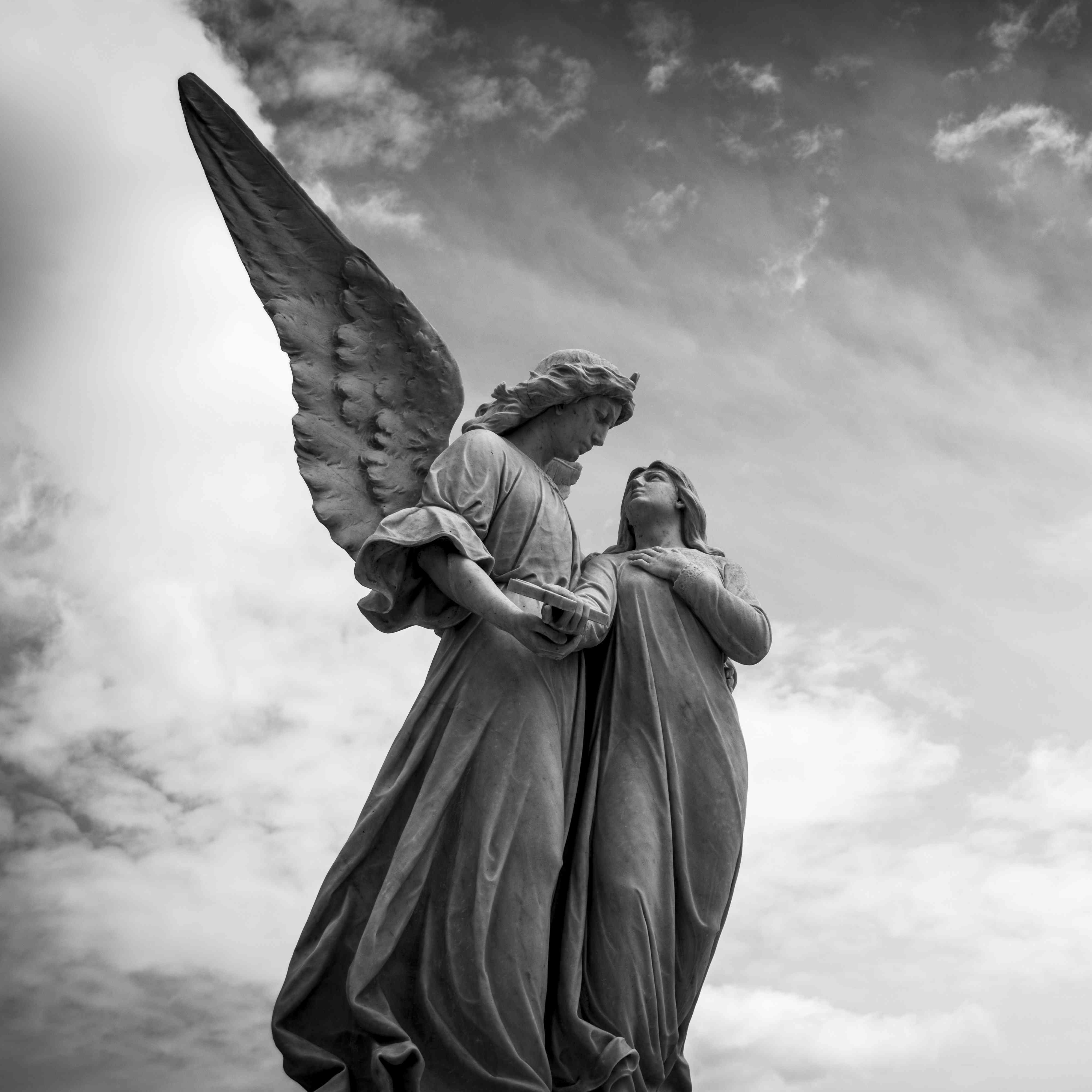}}}\
  \subfloat[Curve length $l$]{{\includegraphics[width=0.23\linewidth]{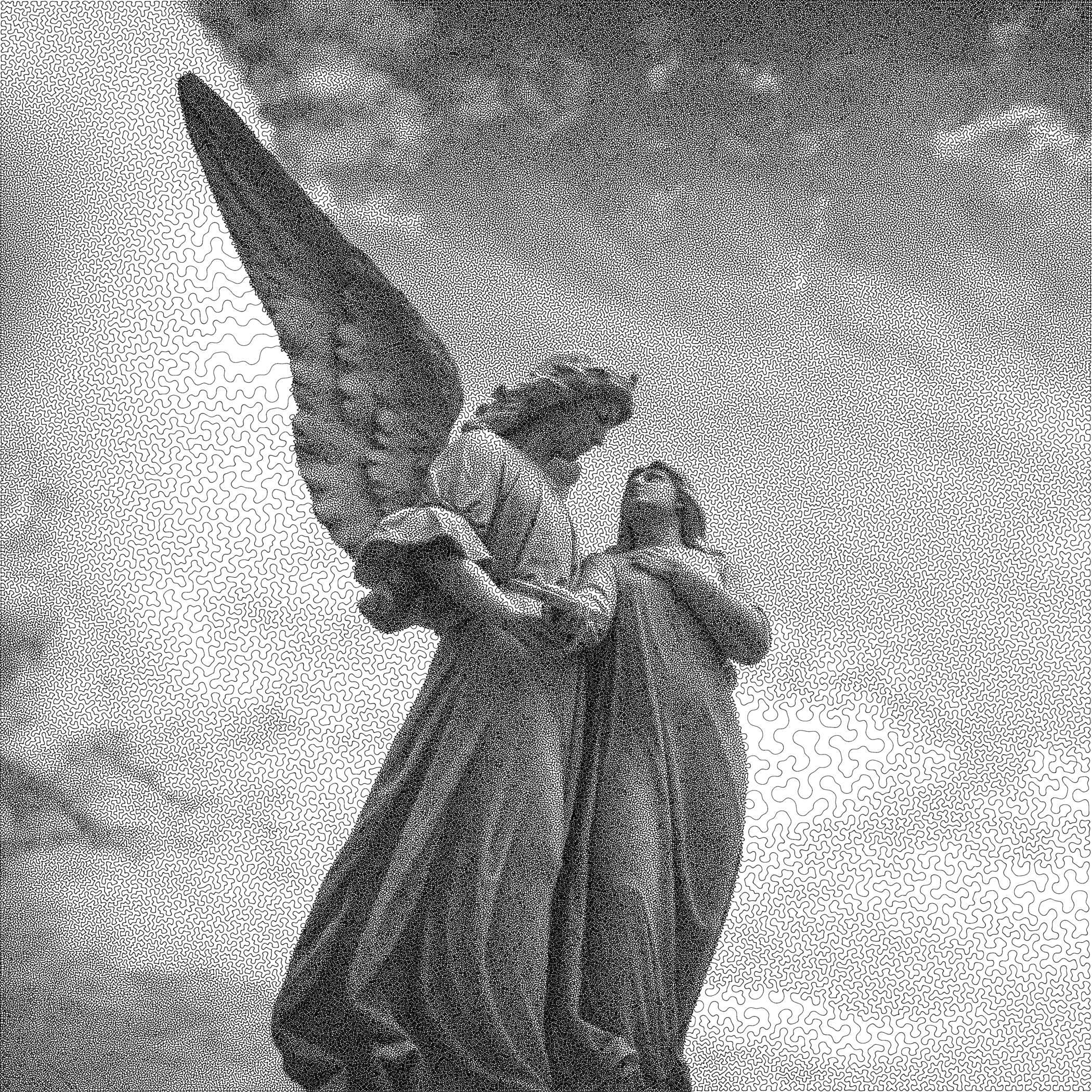}}}\
  \subfloat[Curve length  $\frac{l}{3}$]{{\includegraphics[width=.23\linewidth]{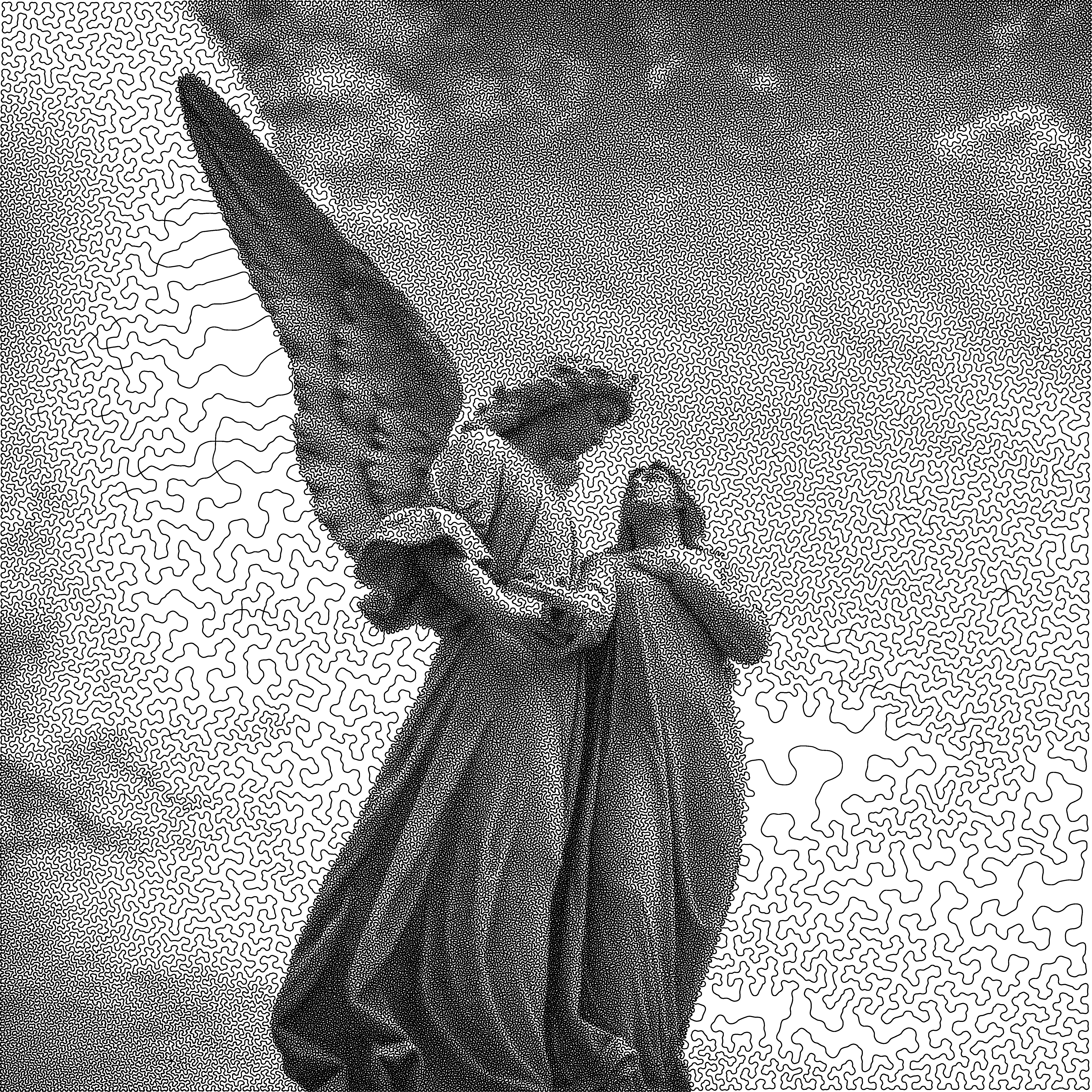}}}\
  \subfloat[Curve length  $\frac{l}{12}$]{{\includegraphics[width=.23\linewidth]{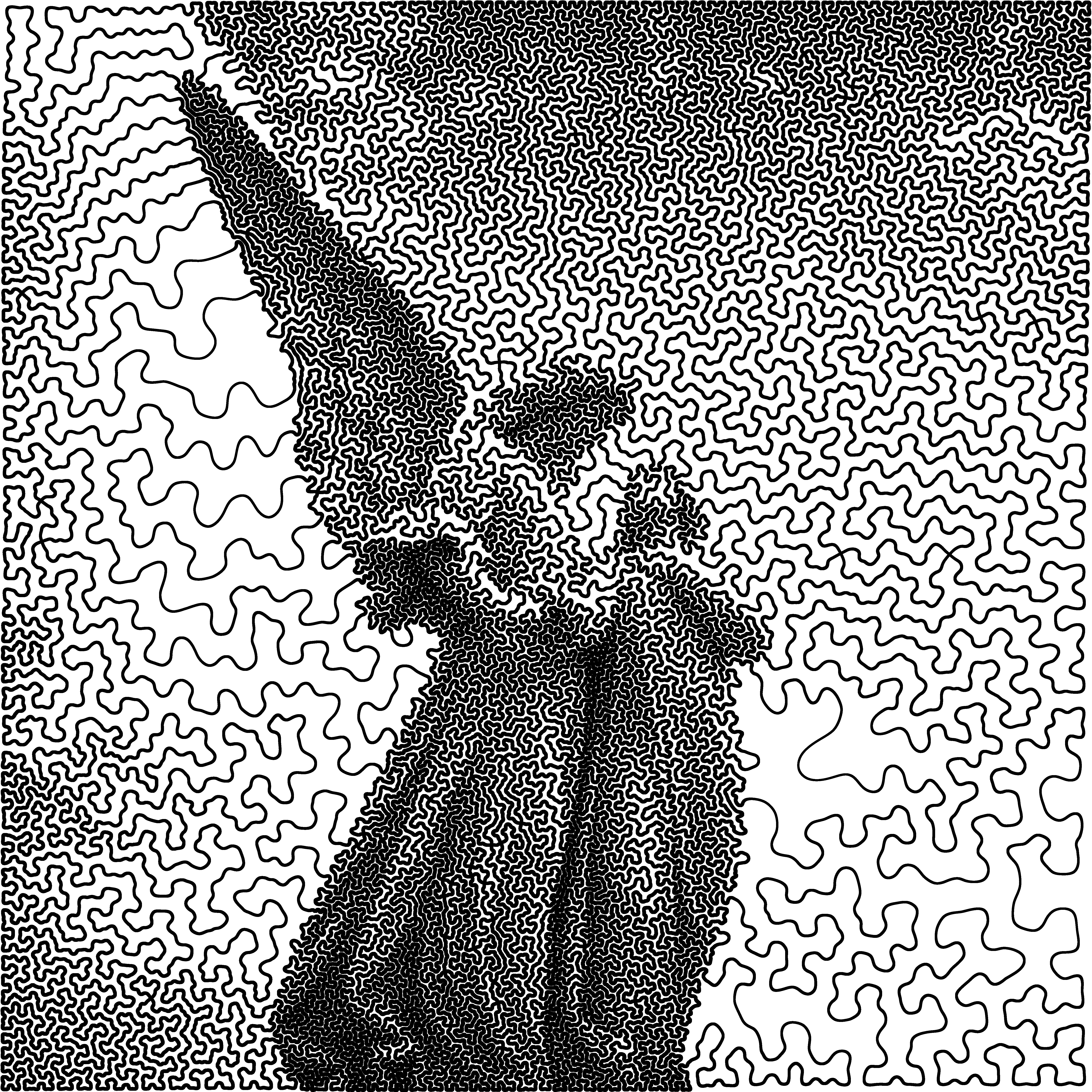}}}
 \caption{Examples of Curvling (stippling + curve projection, 256k, $\approx$ 10'), \label{fig:MadoneConstraint}}
\end{figure*}

A direct application of the proposed algorithm allows to approximate an arbitrary image with measures supported on curves. An example is displayed in Fig. \ref{fig:MadoneConstraint} with curves satisfying different kinematic constraints.

\subsubsection{Color images}

There are different ways to render color images with the proposed idea. We refer for instance to \cite{wei2010multi,chauffert2015comment} for two different examples. 
In this section, we propose a simple alternative idea to give a color to the dots or curves. Given a target vectorial density $\rho=(\rho_R,\rho_G,\rho_B): \Omega \to [0,1]^3$, the algorithm we propose simply reads as follows:
\begin{enumerate}
 \item[1)] We first construct a gray level image defined by:
\begin{equation}
\bar \rho = (\rho_R+\rho_G+\rho_B)/3. 
\end{equation}
 \item[2)] Then, we project the density $\bar \rho$ onto the set of constraints $\mathcal{M}$ with Algorithm \ref{alg:mainalgo}. This yields a sequence of points $\bx\in \Omega^n$.
 \item[3)] Then, for each point $\bx[i]$ of the discretized measure, we select a color as $\frac{\rho(\bx[i])}{\bar \rho(\bx[i]))}$.
\end{enumerate}

\begin{figure*}
  \centering
  \subfloat[Target color image]{{\includegraphics[width=.4\textwidth]{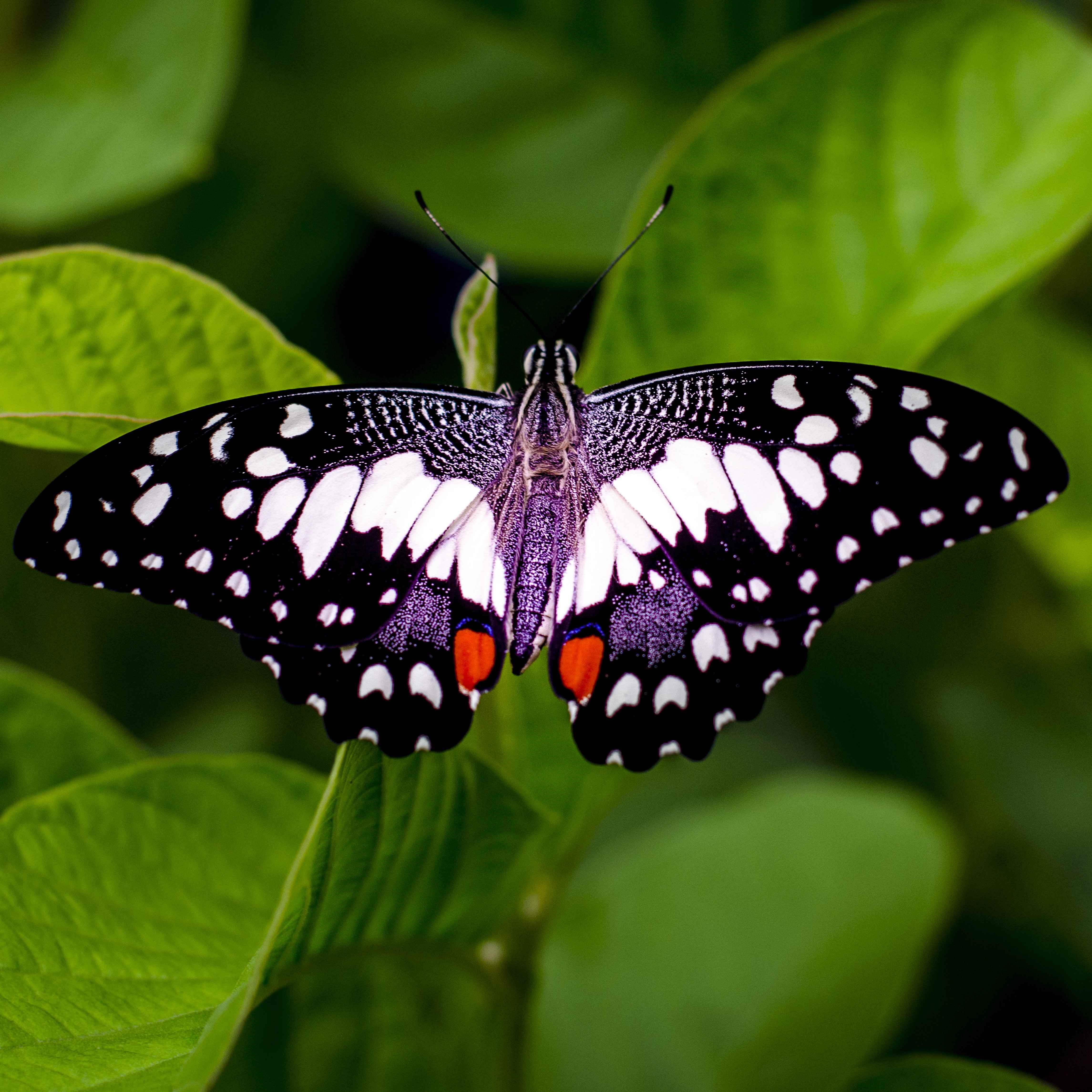}}} \ 
  \subfloat[Approximate color measure]{{\includegraphics[width=.4\textwidth]{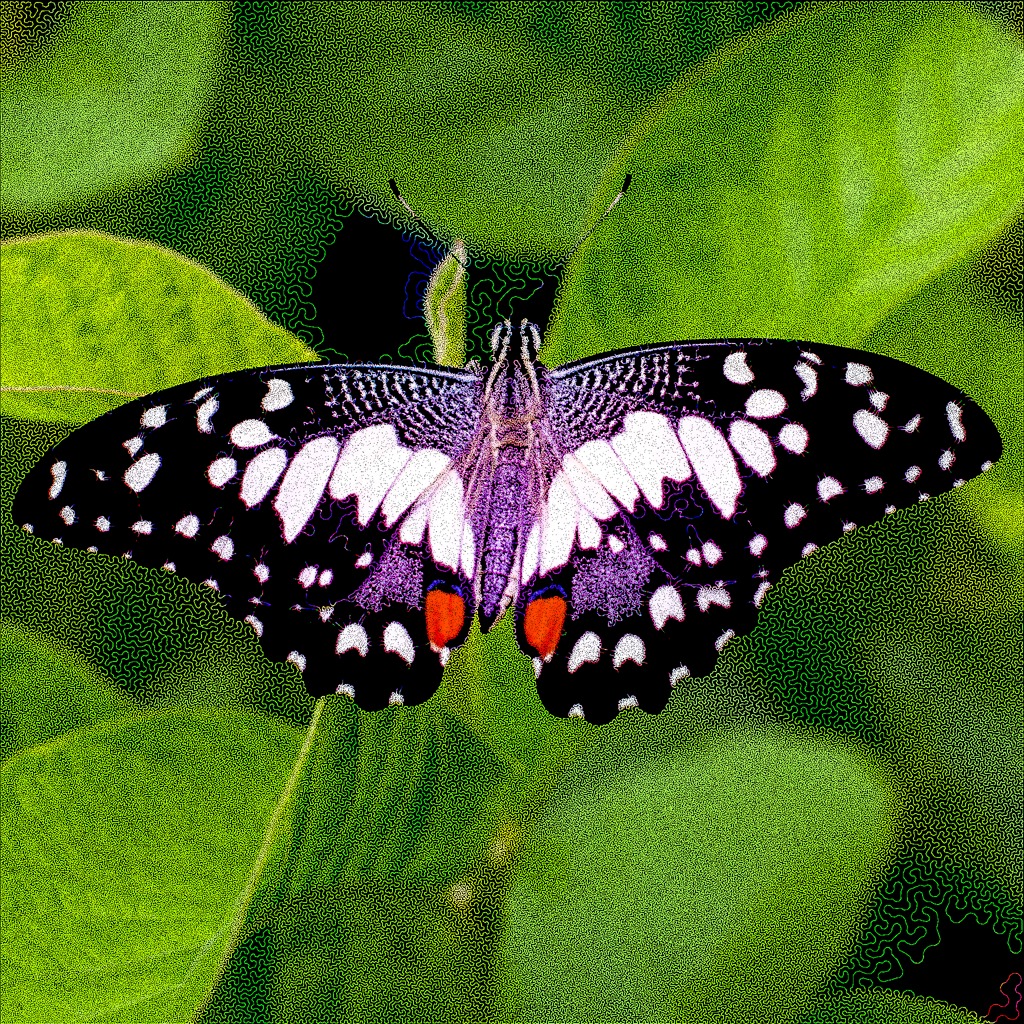}}}
 \caption{Examples of color curvling, 512k, $\approx$ 24'), \label{fig:ColorCurve}}
\end{figure*}

We use only saturated colors, explaining the division in step 3). The parallel for gray-scale images, is that we represent stippling results with disks taking only the maximal intensity. Then, the mean in step 1) is used to attract the curve towards the regions of high luminance of the image.
An example of result of the proposed algorithm is shown in Figure \ref{fig:ColorCurve}.

\subsubsection{Dynamic examples}

The codes can also be used to approximate videos. The principle is simple: first we approximate the first sequence of the frame with our projection algorithm starting from an arbitrary initial guess. Then, the other frames are obtained with the projection algorithm, taking as an initial guess, the result of the previous iteration. This ensures some continuity of the dots or curves between consecutive frames.
Some videos are given in the supplementary material.

\subsection{Path planning}

In this section, we provide two applications of the proposed algorithm to path planning problems. 

\subsubsection{Videodrone}

\begin{figure*}[ht]
  \begin{center}
      \subfloat[The crime density $\mu$]{\includegraphics[width=0.47\textwidth]{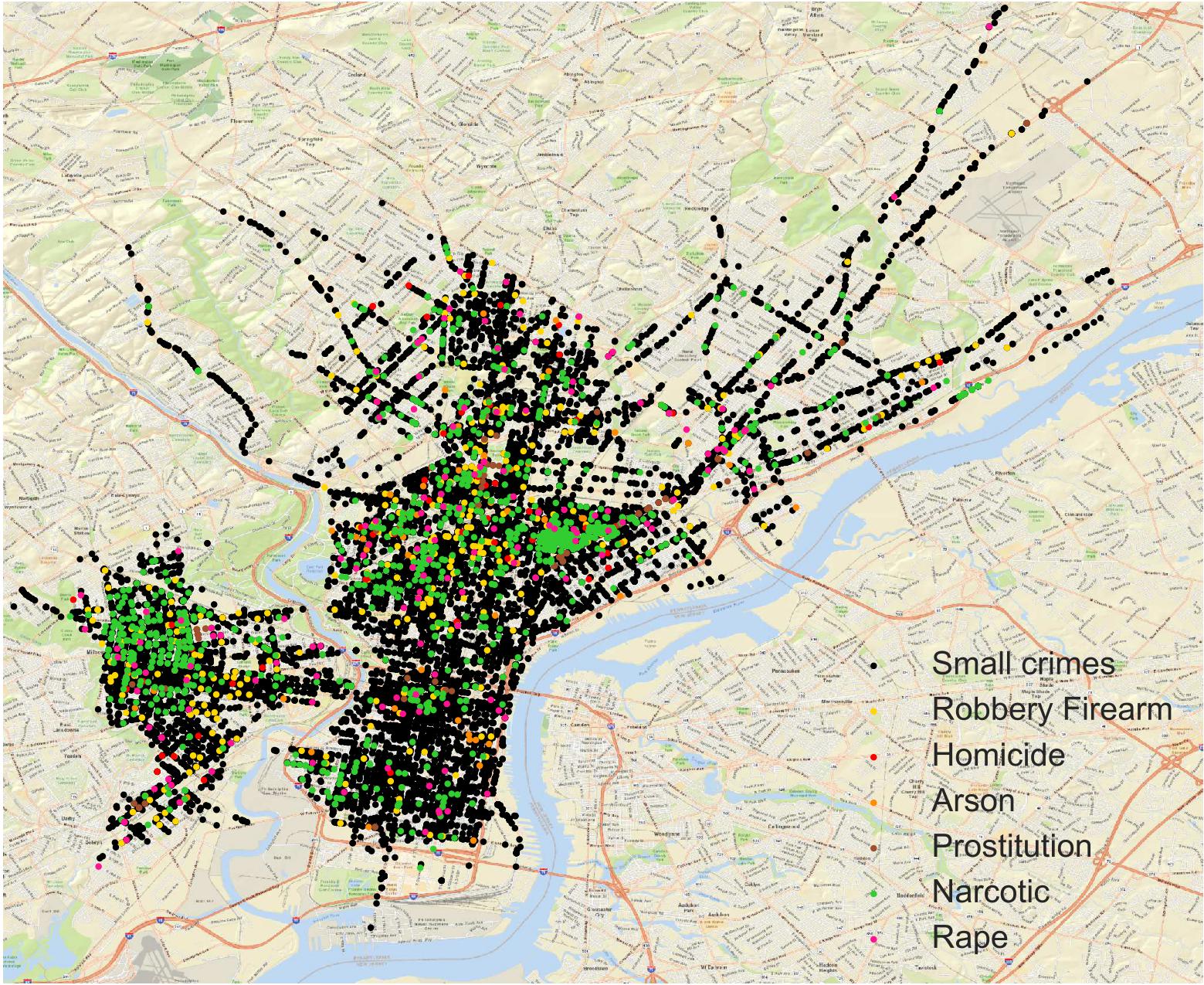}\label{sub:crimeZone}}\ 
      \subfloat[Path adopted by the drone]{\includegraphics[width=0.47\textwidth]{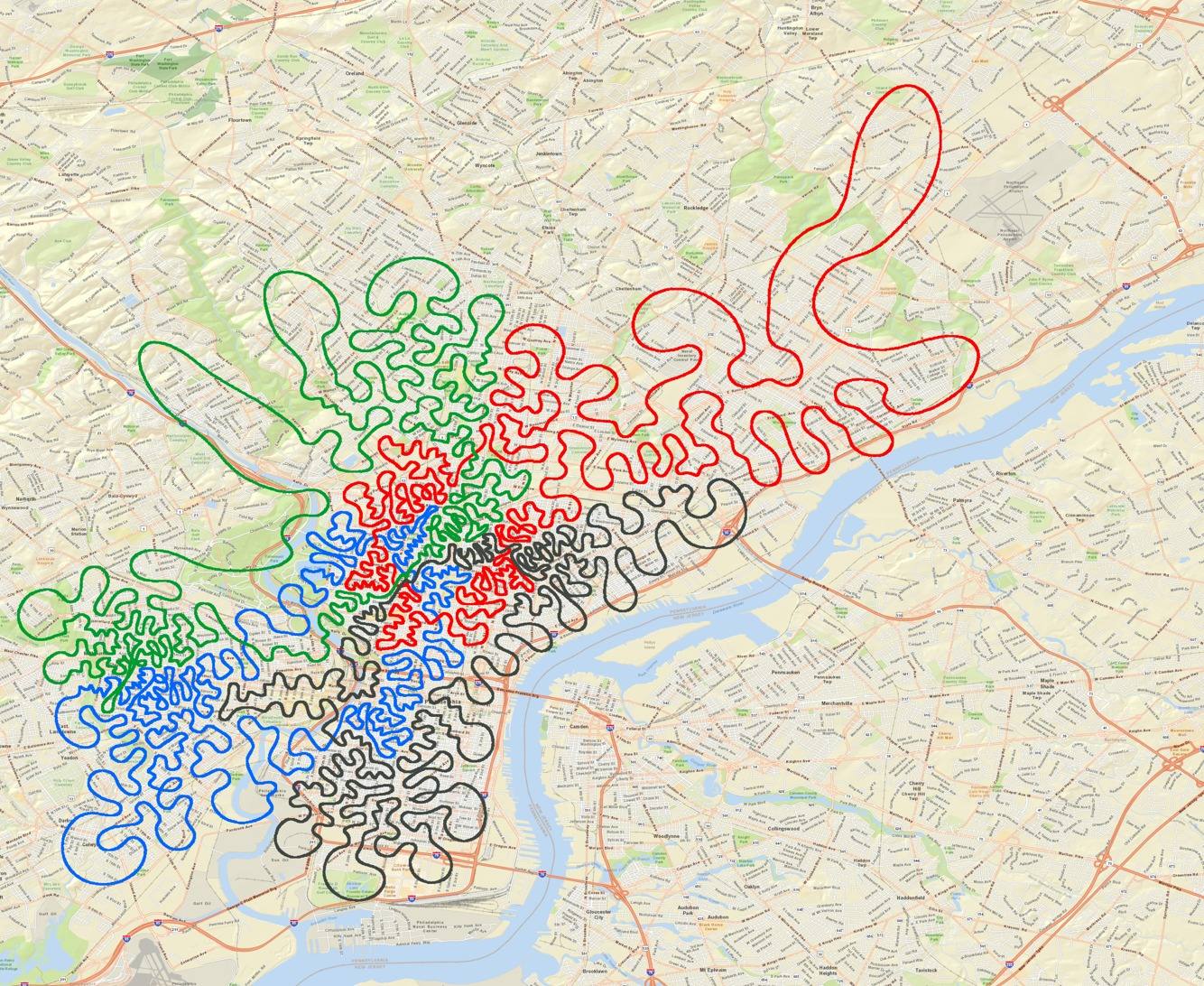}\label{sub:DroneFlight}}
      \caption{The data superimposed on a map of Philadelphia. A possible drone trajectory made. In this example, the drone passes 4 times to its recharging location, explaining the different colors of  the trajectory. In this example, the trajectory was discretized with 8k points and optimized in 30''.}
  \end{center}
\end{figure*}

Drone surveillance is an application with increasing interest from cities, companies or even private individuals.
In this paragraph, we show that the proposed algorithms can be used to plan the drone trajectories for surveillance applications.
We use the criminal data provided by \cite{phily1} to create a density map of crime in Philadelphia, see Fig. \ref{sub:crimeZone}. We give different weights to different types of crimes.  By minimizing \eqref{eq:mainproblem0}, we can design an optimal path, in the sense that it satisfies the kinematic constraints of the drone and passes close to dangerous spots more often than in the remaining locations. 
In this example, we impose a bounded speed, a maximal yaw angular velocity and also to pass at a given location at a given time to recharge the drone to satisfy autonomy constraints.

\subsubsection{Laser engraving}

\begin{figure}[ht]
  \begin{center}
    \subfloat[Laser engraving machine]{\includegraphics[width=0.34\textwidth]{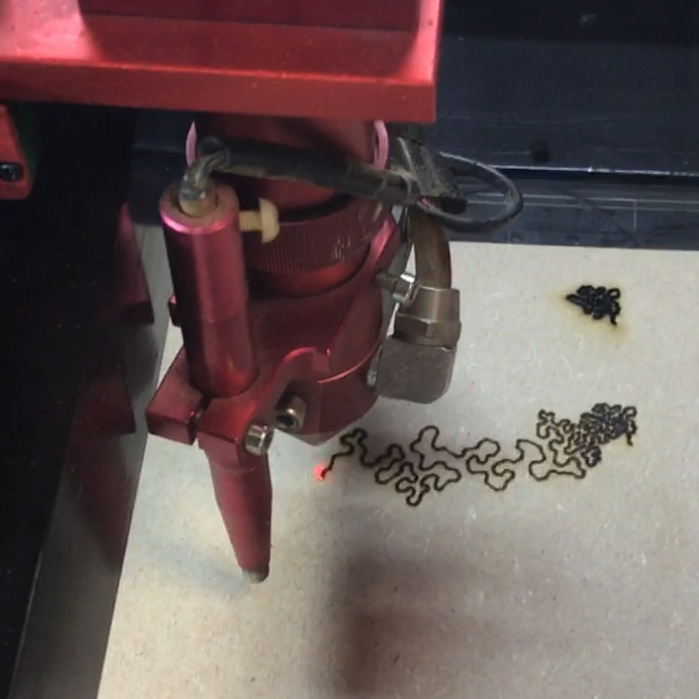}\label{sub:laserengrave1}} 
    \subfloat[The resulting wood engraved trajectory]{\includegraphics[width=0.61\textwidth]{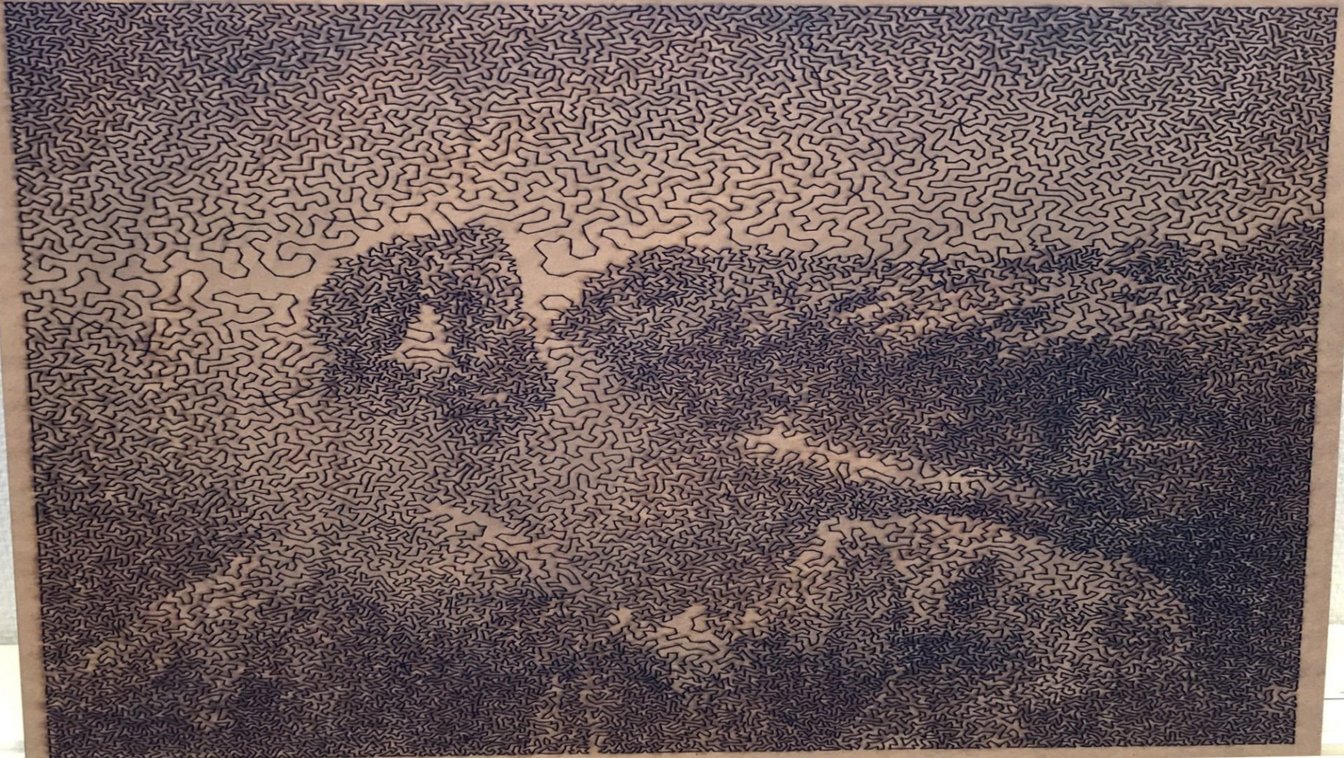}\label{sub:laserengrave2}}\\
      \caption{Example of wood engraving. Left: a laser burning the wood by following an input trajectory. Right: the final result.}
    \label{fig:laser}
  \end{center}
\end{figure}

In Fig. \ref{fig:laser}, we gave a trajectory to a laser engraving machine in order to reproduce a landscape with a continuous line. 
We suspect that the same techniques could be used to optimize the nozzle and the flow of matter trajectory of 3D printers.

\subsubsection{Sampling in MRI}

Following \cite{boyer2016generation}, we propose to generate compressive sampling schemes in MRI (Magnetic Resonance Imaging), using the proposed algorithm.

In MRI, images are probed indirectly through their Fourier transform. 
Fourier transform values are sampled along curves with bounded speed and bounded acceleration, which exactly corresponds to the set of constraints defined in \eqref{eq:kinematicset}. 
The latest compressed sensing theories suggest that a good way of subsampling the Fourier domain, consists in drawing points independently at random according to a certain distribution $\mu$, that depends on the image sparsity structure in the wavelet domain \cite{boyer2016generation,adcock2017breaking}. Unfortunately, this strategy is impractical in MRI due to physical constraints. 
To simulate such a sampling scheme, we therefore propose to project $\mu$ onto the set of admissible trajectories. 

Let $u:[0,1]^2\to \C$ denote a magnetic resonance image. The sampling process yields a set of Fourier transform values $\by[i]= \hat u(\bx[i])$. Given this set of values, the image is then reconstructed by solving a nonlinear convex programming problem:
\begin{equation}
 \min_{v, v|_\bx = \by} \frac{1}{2}\|\hat v(\bx) - \by\|_2^2 + \lambda \|\Psi u\|_1,
\end{equation}
where $\Psi$ is a linear sparsifying transform, such as a redundant wavelet transform. 

\begin{figure}[ht]
  \begin{center}
    \subfloat[Target density $\mu$]{\includegraphics[width=0.33\textwidth]{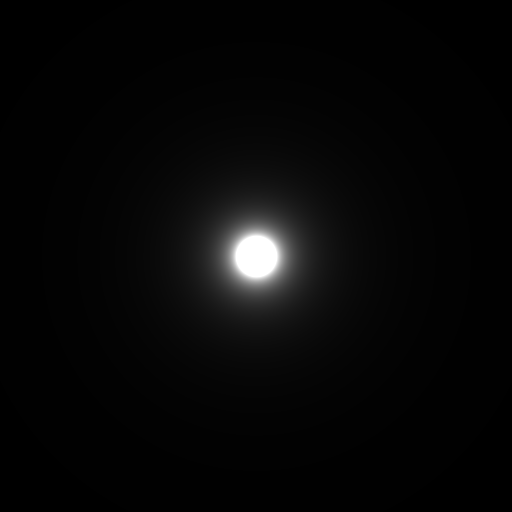}\label{sub:density}} \ 
    \subfloat[Sampling scheme]{\includegraphics[width=0.33\textwidth]{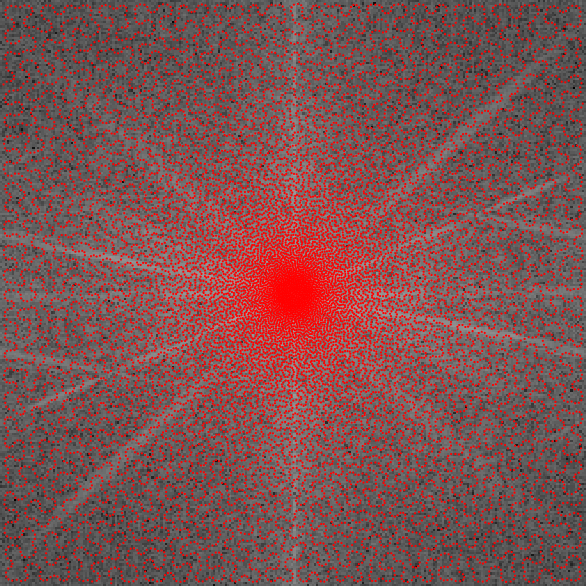}\label{sub:sampling}}\\
    \subfloat[True image]{\includegraphics[width=0.33\textwidth]{MRI_Original.png}\label{sub:original}}\ 
    \subfloat[Reconstructed image]{\includegraphics[width=0.33\textwidth]{MRI_Reconstructed.png}\label{sub:reconstructed}}
      \caption{Example of sampling scheme generation and image reconstruction in MRI. The target density $\mu$ is shown in \ref{sub:density}. The sampling scheme generated by our algorithm is shown in \ref{sub:sampling}. The background shows the Fourier transform of $u$ in log-scale. It contains one fourth of the total number of Fourier transform values. The true image and the reconstructed one are shown in Fig. \ref{sub:original} and \ref{sub:reconstructed}.}
    \label{fig:renonculacees}
  \end{center}
\end{figure}

\clearpage
\appendix
\section{Theoretical convergence of Algorithm \ref{alg:mainalgo}}
\label{sec:app_convergence}

The following result is a direct application of standard convergence results, see e.g. \cite{nesterov2013gradient}.
\begin{theorem}\label{thm:convergencefirstorder}
 Suppose that $X\subset \R^n$ is closed and convex, that $\Sigma_k=\Sigma$ is a constant positive definite matrix. In addition, suppose that $F$ is a $C^1$ function with Lipschitz continuous gradient:
 \begin{equation}\label{eq:defLipschitz}
  \forall (\bx_1,\bx_2), \|\nabla F(\bx_1) - \nabla F(\bx_2)\|_{\Sigma^{-1}} \leq L\|\bx_1-\bx_2\|_{\Sigma}.
 \end{equation}
 Finally suppose that either $X$ is compact or $F$ is coercive. Then Algorithm \ref{alg:mainalgo} converges to a critical point of $F$ for step-size $s_k = \frac 1L$.
\end{theorem}

Applying Theorem \ref{thm:convergencefirstorder} requires  $\Sigma_k$ to be constant, hence the mass $\bw$ to be prescribed. We make this assumption in this section.

Theorem \ref{thm:convergencefirstorder} shows that it is critical to evaluate - if it exists - the Lipschitz constant of $\nabla_\bx F$. By equation \eqref{eq:derivativex}, we need to evaluate the variations of the Laguerre cells barycenter $\bb$ with respect to $\bx$. Unfortunately, following the Hessian computation in \cite{Fred2018differentiation}, the Lipschitz constant scales as $\min_{i\ne j}\Vert \bx[i]-\bx[j]\Vert^{-1}$ and cannot be proven to be uniform in $\bx$. Hence, we can only hope for a local result describing the Lipschitz constant.

Hence, if the hypothesis of existence of the Hessian of $F$ are met (see \cite{Fred2018differentiation}), an estimation of the Lipschitz constant of $F$ by its Hessian yields a theory of local convergence of $F$ in a vicinity $\mathcal{V}(\bx^*)$ of a local minimizer $\bx^*$, for small enough steps $s_k$. Without these assumptions, local Lipschitz continuity of the gradient of $F$ cannot be enforced. 

If in addition there is no $\Pi$-step, that is $X=\Omega^n$, the gradient of $F$ is $1$-Lipschitz around critical points (the so-called centroidal tessellation), see \cite[Prop. 6.3]{du1999centroidal}. Hence, convergence can be proven in a vicinity of $\bx^*$ for step choice $s_k=1$ and the metric $\Sigma_k$. However, the size of the vicinity $\mathcal{V}(\bx^*)$ relies on the geometrical properties of the ``optimal'' Laguerre tessellation. The quality of such a local minimum could be very far from the global minimizer, nevertheless numerical experiments tend to indicate that is it not the case. For the same problem, hundreds of random initializations converge to a set of stationary points with a nice visual property.

\section*{Acknowledgments}
  The authors wish to thank Alban Gossard warmly for his help in designing numerical integration procedures.

\bibliographystyle{abbrv} 

\bibliography{Biblio}   


\end{document}

%% file: 5pts.tex
\begin{tikzpicture}[scale=4.5]
    \coordinate (NW) at (0,1);
    \coordinate (NE) at (1,1);
    \coordinate (SE) at (1,0);
    \coordinate (SW) at (0,0);
    \draw (NW) -- (NE) -- (SE)-- (SW) -- cycle;
     \foreach \i in {0,...,4}
       \draw (.5,\i*.2+.1) circle(0.005);
     \foreach \i in {0,...,5}
       \draw [line width=0.10mm] (0,\i*.2) -- (1,\i*.2);
       \draw [ultra thin,<->] (.5,.5) -- (.5,.70) node[pos=0.65,right,scale=0.8]{$\| x_i - x_{i+1}\|$};
\end{tikzpicture}

%% file: 25pts.tex
\begin{tikzpicture}[scale=4.5]
    \coordinate (NW) at (0,1);
    \coordinate (NE) at (1,1);
    \coordinate (SE) at (1,0);
    \coordinate (SW) at (0,0);
    \draw (NW) -- (NE) -- (SE)-- (SW) -- cycle;
     \foreach \i in {0,...,24}
       \draw (.5,\i*.04+.02) circle(0.005) ;
     \foreach \i in {0,...,25}
       \draw [line width=0.10mm] (0,\i*.04) -- (1,\i*.04);
      
     \draw [ultra thin,<->] (.5,.5) -- (.5,.54) node[pos=0.9,right,scale=0.8]{$\| x_i - x_{i+1}\|$};
\end{tikzpicture}

%% file: Cat_ini.tex
%
%
\begin{tikzpicture}

\begin{axis}[%
width=0.3\textwidth,
height=0.3\textwidth,
at={(0.766in,0.486in)},
scale only axis,
xmin=31,
xmax=570,
ymin=-18.0572580645161,
ymax=407.057258064516,
axis lines=none, xtick=\empty, ytick=\empty,
axis background/.style={fill=white}
]
\addplot [color=red,solid,forget plot]
  table[row sep=crcr]{%
31	384\\
32	385\\
32	386\\
33	387\\
34	387\\
35	388\\
36	388\\
37	389\\
38	389\\
39	389\\
40	389\\
41	389\\
42	389\\
43	389\\
44	389\\
45	389\\
46	389\\
47	389\\
48	389\\
49	389\\
50	389\\
51	388\\
52	388\\
53	388\\
54	388\\
55	388\\
56	388\\
57	387\\
58	387\\
59	387\\
60	387\\
61	386\\
62	386\\
63	385\\
64	384\\
65	384\\
66	383\\
67	383\\
68	382\\
69	381\\
70	381\\
71	380\\
72	379\\
73	378\\
74	377\\
74	376\\
75	375\\
76	374\\
76	373\\
77	372\\
78	371\\
78	370\\
78	369\\
79	368\\
79	367\\
80	366\\
81	365\\
81	364\\
82	363\\
82	362\\
83	361\\
83	360\\
84	359\\
84	358\\
84	357\\
85	356\\
85	355\\
86	354\\
86	353\\
87	352\\
87	351\\
88	350\\
88	349\\
88	348\\
89	347\\
89	346\\
90	345\\
90	344\\
91	343\\
91	342\\
92	341\\
92	340\\
93	339\\
93	338\\
94	337\\
95	336\\
95	335\\
96	334\\
96	333\\
97	332\\
97	331\\
98	330\\
98	329\\
99	328\\
100	327\\
100	326\\
101	325\\
101	324\\
102	323\\
103	322\\
103	321\\
104	320\\
104	319\\
105	318\\
105	317\\
106	316\\
107	315\\
107	314\\
108	313\\
109	312\\
110	311\\
111	310\\
111	309\\
112	308\\
113	307\\
114	306\\
114	305\\
115	304\\
116	303\\
117	302\\
117	301\\
118	300\\
119	299\\
120	298\\
120	297\\
121	296\\
122	295\\
123	294\\
124	293\\
124	292\\
125	291\\
126	290\\
127	289\\
128	288\\
129	287\\
130	286\\
131	285\\
132	284\\
133	283\\
134	282\\
135	281\\
136	280\\
137	279\\
138	278\\
139	278\\
140	277\\
140	276\\
141	275\\
142	274\\
143	273\\
144	273\\
145	272\\
146	271\\
147	270\\
148	269\\
149	268\\
150	268\\
151	267\\
152	266\\
153	265\\
154	264\\
155	264\\
156	263\\
157	262\\
158	261\\
159	261\\
160	260\\
161	260\\
162	259\\
163	259\\
164	259\\
165	258\\
166	258\\
167	258\\
168	257\\
169	257\\
170	257\\
171	256\\
172	256\\
173	255\\
174	255\\
175	255\\
176	254\\
177	254\\
178	254\\
179	253\\
180	254\\
180	255\\
181	255\\
182	255\\
183	254\\
184	254\\
185	254\\
186	254\\
187	253\\
188	253\\
189	253\\
190	253\\
191	252\\
192	252\\
193	252\\
194	252\\
195	252\\
196	252\\
197	252\\
198	253\\
199	253\\
200	253\\
201	253\\
202	253\\
203	253\\
204	253\\
205	253\\
206	253\\
207	253\\
208	253\\
209	253\\
210	253\\
211	253\\
212	253\\
213	253\\
214	253\\
215	253\\
216	253\\
217	253\\
218	253\\
219	253\\
220	253\\
221	253\\
222	253\\
223	253\\
224	253\\
225	253\\
226	253\\
227	253\\
228	253\\
229	253\\
230	253\\
231	252\\
232	252\\
233	252\\
234	252\\
235	252\\
236	252\\
237	252\\
238	252\\
239	252\\
240	252\\
241	251\\
242	251\\
243	251\\
244	251\\
245	251\\
246	251\\
247	250\\
248	250\\
249	250\\
250	250\\
251	250\\
252	250\\
253	250\\
254	249\\
255	249\\
256	249\\
257	249\\
258	249\\
259	249\\
260	248\\
261	248\\
262	248\\
263	248\\
264	248\\
265	248\\
266	247\\
267	247\\
268	247\\
269	247\\
270	247\\
271	247\\
272	246\\
273	246\\
274	246\\
275	246\\
276	246\\
277	246\\
278	246\\
279	246\\
280	246\\
281	245\\
282	245\\
283	245\\
284	245\\
285	245\\
286	245\\
287	245\\
288	245\\
289	244\\
290	244\\
291	244\\
292	244\\
293	244\\
294	244\\
295	244\\
296	244\\
297	244\\
298	243\\
299	243\\
300	243\\
301	243\\
302	243\\
303	243\\
304	243\\
305	243\\
306	242\\
307	242\\
308	242\\
309	242\\
310	242\\
311	242\\
312	242\\
313	241\\
314	241\\
315	241\\
316	241\\
317	241\\
318	241\\
319	241\\
320	240\\
321	240\\
322	240\\
323	240\\
324	240\\
325	240\\
326	240\\
327	240\\
328	239\\
329	239\\
330	239\\
331	239\\
332	239\\
333	239\\
334	239\\
335	239\\
336	239\\
337	239\\
338	239\\
339	239\\
340	238\\
341	238\\
342	238\\
343	238\\
344	238\\
345	238\\
346	238\\
347	238\\
348	238\\
349	238\\
350	238\\
351	238\\
352	238\\
353	238\\
354	238\\
355	238\\
356	238\\
357	238\\
358	238\\
359	238\\
360	238\\
361	238\\
362	239\\
363	239\\
364	239\\
365	239\\
366	239\\
367	239\\
368	239\\
369	239\\
370	239\\
371	239\\
372	239\\
373	239\\
374	239\\
375	240\\
376	240\\
377	240\\
378	240\\
379	240\\
380	240\\
381	240\\
382	240\\
383	241\\
384	241\\
385	241\\
386	241\\
387	241\\
388	242\\
389	242\\
390	242\\
391	242\\
392	242\\
393	242\\
394	243\\
395	243\\
396	243\\
397	244\\
398	244\\
399	244\\
400	245\\
401	245\\
402	245\\
403	246\\
404	246\\
405	246\\
406	247\\
407	247\\
408	247\\
409	247\\
410	247\\
411	248\\
412	248\\
413	248\\
414	248\\
415	248\\
416	248\\
417	248\\
418	248\\
419	248\\
420	249\\
421	249\\
422	249\\
423	249\\
424	249\\
425	249\\
426	249\\
427	249\\
428	249\\
429	249\\
430	249\\
431	248\\
432	248\\
433	248\\
434	248\\
435	248\\
436	248\\
437	248\\
438	248\\
439	248\\
440	247\\
441	247\\
442	247\\
443	247\\
444	247\\
445	247\\
446	247\\
447	247\\
448	247\\
449	246\\
450	246\\
451	246\\
452	246\\
453	246\\
454	246\\
455	247\\
456	247\\
457	247\\
458	247\\
459	247\\
460	247\\
461	246\\
462	246\\
463	246\\
464	246\\
465	246\\
466	247\\
467	247\\
468	247\\
469	247\\
470	248\\
471	248\\
472	248\\
473	249\\
474	249\\
475	249\\
476	249\\
477	250\\
478	250\\
479	250\\
480	250\\
481	251\\
482	251\\
483	251\\
484	252\\
485	253\\
486	254\\
487	254\\
488	255\\
489	256\\
490	257\\
491	257\\
492	258\\
493	259\\
494	260\\
495	260\\
496	261\\
497	262\\
498	262\\
499	263\\
500	264\\
501	264\\
502	264\\
503	265\\
504	265\\
505	266\\
506	266\\
507	267\\
508	267\\
509	268\\
510	269\\
510	270\\
511	271\\
512	272\\
512	273\\
513	274\\
514	275\\
515	276\\
516	276\\
517	277\\
518	278\\
519	279\\
520	280\\
521	281\\
522	282\\
523	283\\
523	284\\
524	285\\
525	286\\
525	287\\
526	288\\
527	289\\
528	290\\
529	290\\
530	290\\
531	289\\
531	288\\
531	287\\
531	286\\
531	285\\
531	284\\
530	283\\
530	282\\
530	281\\
530	280\\
531	279\\
532	280\\
533	280\\
534	281\\
535	282\\
536	283\\
537	283\\
538	284\\
539	285\\
540	286\\
541	287\\
542	287\\
543	288\\
544	289\\
545	289\\
545	288\\
545	287\\
545	286\\
545	285\\
545	284\\
545	283\\
545	282\\
544	281\\
544	280\\
544	279\\
544	278\\
544	277\\
544	276\\
543	275\\
543	274\\
543	273\\
542	272\\
542	271\\
542	270\\
542	269\\
541	268\\
541	267\\
540	266\\
540	265\\
540	264\\
539	263\\
539	262\\
539	261\\
540	260\\
541	260\\
542	259\\
543	258\\
544	258\\
545	257\\
546	256\\
547	256\\
547	255\\
548	254\\
549	253\\
550	252\\
551	251\\
552	250\\
553	249\\
553	248\\
554	247\\
555	246\\
556	245\\
557	244\\
557	243\\
558	242\\
558	241\\
558	240\\
558	239\\
558	238\\
558	237\\
559	236\\
560	235\\
560	234\\
561	233\\
561	232\\
562	231\\
563	230\\
563	229\\
564	228\\
564	227\\
565	226\\
565	225\\
566	224\\
566	223\\
567	222\\
567	221\\
568	220\\
568	219\\
568	218\\
569	217\\
569	216\\
569	215\\
570	214\\
570	213\\
570	212\\
569	211\\
569	210\\
568	209\\
567	208\\
567	207\\
566	206\\
566	205\\
566	204\\
565	203\\
565	202\\
564	201\\
563	200\\
563	199\\
562	199\\
561	198\\
560	197\\
560	196\\
559	196\\
558	196\\
557	195\\
556	195\\
555	195\\
554	195\\
553	194\\
552	193\\
551	193\\
550	192\\
549	192\\
548	192\\
547	192\\
546	192\\
545	192\\
544	192\\
543	192\\
542	193\\
541	193\\
540	193\\
539	193\\
538	193\\
537	193\\
536	194\\
535	194\\
534	194\\
533	194\\
532	194\\
531	194\\
530	194\\
529	194\\
528	194\\
527	194\\
526	194\\
525	194\\
524	194\\
523	194\\
522	194\\
521	194\\
520	194\\
519	194\\
518	194\\
517	194\\
516	193\\
515	193\\
514	193\\
513	193\\
512	192\\
511	191\\
510	191\\
509	190\\
508	189\\
507	189\\
506	188\\
505	188\\
504	187\\
504	188\\
503	189\\
502	188\\
501	187\\
500	186\\
499	185\\
498	184\\
498	183\\
497	182\\
496	181\\
495	180\\
494	179\\
493	180\\
492	179\\
492	178\\
491	177\\
491	176\\
491	175\\
490	174\\
490	173\\
490	172\\
489	171\\
489	170\\
489	169\\
488	168\\
488	167\\
487	166\\
486	165\\
486	164\\
485	163\\
484	162\\
483	161\\
483	160\\
482	159\\
481	158\\
480	157\\
479	156\\
478	155\\
477	154\\
477	153\\
476	152\\
476	151\\
475	150\\
475	149\\
474	148\\
474	147\\
473	146\\
473	145\\
472	144\\
472	143\\
471	142\\
470	142\\
469	143\\
468	142\\
468	141\\
467	140\\
467	139\\
467	138\\
467	137\\
466	136\\
466	135\\
466	134\\
465	133\\
465	132\\
464	131\\
464	130\\
463	129\\
463	128\\
462	127\\
461	127\\
460	126\\
460	125\\
460	124\\
460	123\\
460	122\\
460	121\\
460	120\\
460	119\\
460	118\\
460	117\\
461	116\\
461	115\\
462	114\\
462	113\\
463	112\\
464	111\\
464	110\\
465	109\\
466	108\\
467	107\\
467	106\\
468	105\\
469	104\\
470	103\\
471	102\\
472	101\\
473	100\\
474	99\\
475	98\\
476	97\\
477	96\\
478	95\\
479	94\\
480	94\\
481	93\\
482	92\\
482	91\\
483	90\\
484	89\\
485	88\\
486	87\\
487	86\\
488	85\\
489	84\\
490	83\\
490	82\\
491	81\\
492	80\\
492	79\\
493	78\\
494	77\\
494	76\\
494	75\\
495	74\\
495	73\\
496	72\\
496	71\\
497	70\\
497	69\\
498	68\\
498	67\\
499	66\\
499	65\\
500	64\\
500	63\\
501	62\\
501	61\\
502	60\\
502	59\\
503	58\\
504	57\\
504	56\\
505	55\\
505	54\\
506	53\\
506	52\\
507	51\\
507	50\\
508	49\\
508	48\\
509	47\\
510	46\\
510	45\\
510	44\\
511	43\\
511	42\\
511	41\\
512	40\\
512	39\\
512	38\\
513	37\\
513	36\\
513	35\\
514	34\\
514	33\\
514	32\\
515	31\\
515	30\\
515	29\\
516	28\\
516	27\\
516	26\\
516	25\\
516	24\\
516	23\\
516	22\\
516	21\\
516	20\\
516	19\\
516	18\\
516	17\\
516	16\\
516	15\\
516	14\\
516	13\\
516	12\\
516	11\\
516	10\\
515	9\\
515	8\\
515	7\\
514	6\\
513	6\\
512	5\\
511	5\\
510	6\\
509	6\\
508	6\\
507	6\\
506	7\\
505	7\\
504	8\\
503	9\\
502	10\\
501	11\\
500	12\\
499	13\\
498	14\\
497	15\\
496	15\\
495	16\\
494	17\\
494	18\\
493	19\\
493	20\\
493	21\\
492	22\\
492	23\\
492	24\\
491	25\\
491	26\\
490	27\\
490	28\\
490	29\\
489	30\\
489	31\\
489	32\\
488	33\\
488	34\\
487	35\\
487	36\\
487	37\\
486	38\\
486	39\\
486	40\\
485	41\\
485	42\\
484	43\\
484	44\\
484	45\\
483	46\\
483	47\\
482	48\\
482	49\\
482	50\\
481	51\\
481	52\\
480	53\\
479	54\\
479	55\\
478	56\\
477	57\\
477	58\\
476	59\\
475	60\\
475	61\\
474	62\\
473	63\\
473	64\\
472	65\\
471	66\\
470	67\\
469	68\\
468	69\\
467	69\\
466	70\\
465	70\\
464	71\\
463	71\\
462	72\\
461	73\\
460	73\\
459	74\\
458	74\\
457	75\\
456	75\\
455	76\\
454	76\\
453	77\\
452	78\\
451	78\\
450	79\\
449	79\\
448	80\\
447	80\\
446	81\\
445	81\\
444	82\\
443	83\\
442	83\\
441	84\\
440	84\\
439	85\\
438	85\\
437	86\\
436	87\\
435	87\\
434	88\\
433	88\\
432	89\\
431	89\\
430	90\\
429	90\\
428	91\\
427	92\\
426	92\\
425	93\\
424	94\\
423	95\\
422	96\\
421	96\\
420	95\\
420	94\\
420	93\\
420	92\\
420	91\\
420	90\\
420	89\\
419	88\\
419	87\\
419	86\\
419	85\\
419	84\\
419	83\\
419	82\\
419	81\\
419	80\\
419	79\\
419	78\\
419	77\\
419	76\\
419	75\\
418	74\\
418	73\\
418	72\\
418	71\\
418	70\\
418	69\\
418	68\\
418	67\\
418	66\\
418	65\\
418	64\\
418	63\\
418	62\\
418	61\\
418	60\\
418	59\\
418	58\\
418	57\\
418	56\\
418	55\\
418	54\\
418	53\\
418	52\\
418	51\\
418	50\\
418	49\\
418	48\\
418	47\\
418	46\\
418	45\\
418	44\\
418	43\\
419	42\\
419	41\\
419	40\\
419	39\\
420	38\\
420	37\\
420	36\\
421	35\\
421	34\\
422	33\\
422	32\\
423	31\\
424	30\\
425	29\\
426	28\\
427	27\\
428	27\\
429	26\\
430	25\\
431	24\\
432	24\\
433	23\\
434	22\\
435	22\\
436	21\\
437	21\\
438	20\\
439	19\\
440	19\\
441	18\\
442	17\\
443	16\\
444	16\\
444	15\\
445	14\\
446	13\\
446	12\\
446	11\\
446	10\\
445	9\\
444	8\\
443	8\\
442	8\\
441	8\\
440	7\\
439	7\\
438	7\\
437	7\\
436	7\\
435	7\\
434	7\\
433	7\\
432	7\\
431	7\\
430	7\\
429	7\\
428	7\\
427	7\\
426	7\\
425	7\\
424	7\\
423	7\\
422	7\\
421	8\\
420	8\\
419	8\\
418	8\\
417	8\\
416	9\\
415	9\\
414	9\\
413	10\\
412	10\\
411	11\\
410	12\\
409	13\\
408	14\\
408	15\\
407	16\\
407	17\\
406	18\\
405	19\\
405	20\\
405	21\\
405	22\\
405	23\\
405	24\\
405	25\\
404	26\\
404	27\\
404	28\\
404	29\\
404	30\\
404	31\\
403	32\\
402	33\\
401	34\\
401	35\\
400	36\\
400	37\\
400	38\\
400	39\\
400	40\\
399	41\\
398	41\\
397	40\\
397	41\\
397	42\\
397	43\\
397	44\\
397	45\\
397	46\\
397	47\\
397	48\\
397	49\\
397	50\\
397	51\\
397	52\\
397	53\\
397	54\\
397	55\\
397	56\\
397	57\\
397	58\\
397	59\\
396	60\\
396	61\\
396	62\\
396	63\\
396	64\\
396	65\\
395	66\\
395	67\\
395	68\\
395	69\\
395	70\\
395	71\\
394	72\\
394	73\\
394	74\\
394	75\\
394	76\\
393	77\\
393	78\\
393	79\\
393	80\\
392	81\\
392	82\\
392	83\\
392	84\\
391	85\\
391	86\\
391	87\\
391	88\\
390	89\\
390	90\\
390	91\\
389	92\\
389	93\\
389	94\\
389	95\\
388	96\\
388	97\\
388	98\\
387	99\\
387	100\\
387	101\\
386	102\\
386	103\\
385	104\\
385	105\\
385	106\\
384	107\\
384	108\\
384	109\\
383	110\\
383	111\\
383	112\\
382	113\\
382	114\\
383	115\\
383	116\\
383	117\\
382	118\\
381	119\\
380	120\\
380	121\\
379	122\\
378	122\\
377	122\\
376	123\\
375	123\\
374	123\\
373	123\\
372	123\\
371	123\\
370	123\\
369	123\\
368	123\\
367	123\\
366	122\\
365	122\\
364	122\\
363	122\\
362	122\\
361	122\\
360	122\\
359	121\\
358	121\\
357	121\\
356	121\\
356	122\\
355	123\\
354	122\\
353	122\\
352	122\\
351	122\\
350	122\\
349	121\\
348	121\\
347	121\\
346	120\\
345	120\\
345	121\\
344	122\\
343	122\\
342	122\\
341	122\\
340	122\\
339	122\\
338	122\\
337	122\\
336	122\\
335	122\\
334	122\\
333	122\\
332	122\\
331	121\\
330	121\\
329	120\\
328	121\\
327	121\\
326	121\\
325	120\\
324	120\\
323	119\\
322	119\\
321	119\\
320	118\\
319	118\\
318	118\\
317	117\\
316	117\\
315	116\\
314	116\\
313	115\\
312	114\\
311	114\\
310	113\\
309	113\\
308	112\\
307	111\\
306	111\\
306	112\\
306	113\\
306	114\\
305	115\\
305	116\\
304	117\\
303	117\\
302	117\\
301	116\\
300	116\\
299	116\\
298	116\\
297	116\\
296	116\\
295	116\\
294	116\\
293	116\\
292	116\\
291	116\\
290	116\\
289	117\\
288	117\\
287	117\\
286	117\\
285	117\\
284	117\\
283	117\\
282	117\\
281	117\\
280	117\\
279	116\\
278	116\\
277	116\\
278	117\\
279	118\\
279	119\\
278	120\\
277	120\\
276	119\\
275	119\\
274	119\\
273	119\\
272	119\\
271	119\\
270	119\\
269	120\\
268	121\\
267	121\\
266	121\\
265	122\\
264	122\\
263	122\\
262	121\\
261	121\\
260	121\\
259	120\\
259	119\\
258	118\\
257	117\\
257	116\\
256	115\\
255	114\\
254	113\\
254	112\\
253	111\\
252	110\\
252	109\\
251	108\\
250	107\\
249	106\\
248	105\\
248	104\\
247	103\\
246	102\\
245	101\\
244	100\\
243	99\\
242	98\\
241	97\\
240	96\\
239	95\\
238	94\\
237	93\\
236	92\\
235	91\\
235	90\\
234	89\\
234	88\\
233	87\\
233	86\\
232	85\\
232	84\\
231	83\\
230	82\\
230	81\\
229	80\\
229	79\\
228	78\\
228	77\\
227	76\\
226	75\\
226	74\\
225	73\\
225	72\\
224	71\\
224	70\\
223	69\\
223	68\\
222	67\\
221	66\\
221	65\\
220	64\\
220	63\\
219	62\\
219	61\\
219	60\\
219	59\\
218	58\\
218	57\\
218	56\\
219	55\\
219	54\\
219	53\\
219	52\\
219	51\\
219	50\\
219	49\\
220	48\\
221	47\\
221	46\\
222	45\\
223	44\\
223	43\\
224	42\\
225	41\\
226	40\\
227	39\\
228	38\\
229	37\\
229	36\\
230	35\\
231	34\\
232	34\\
232	33\\
233	32\\
234	31\\
235	31\\
235	30\\
236	29\\
237	28\\
238	27\\
239	26\\
240	25\\
241	24\\
242	23\\
243	23\\
244	22\\
245	21\\
246	21\\
247	21\\
248	20\\
249	20\\
250	20\\
251	20\\
252	20\\
253	19\\
254	19\\
255	19\\
256	19\\
257	18\\
258	18\\
259	17\\
260	17\\
261	16\\
262	16\\
263	15\\
264	15\\
265	14\\
266	13\\
267	12\\
268	11\\
268	10\\
268	9\\
268	8\\
269	7\\
269	6\\
269	5\\
270	4\\
270	3\\
269	3\\
268	3\\
267	3\\
266	2\\
265	2\\
264	2\\
263	2\\
262	2\\
261	1\\
260	1\\
259	1\\
258	1\\
257	1\\
256	1\\
255	1\\
254	1\\
253	1\\
252	1\\
251	1\\
250	1\\
249	1\\
248	0\\
247	0\\
246	0\\
245	0\\
244	0\\
243	1\\
242	1\\
241	1\\
240	1\\
239	1\\
238	2\\
237	2\\
236	2\\
235	2\\
234	2\\
233	2\\
232	3\\
231	4\\
230	4\\
229	5\\
228	5\\
227	6\\
226	6\\
225	7\\
224	7\\
223	8\\
222	9\\
221	10\\
220	11\\
219	12\\
218	13\\
217	14\\
216	15\\
215	16\\
214	17\\
213	18\\
212	19\\
211	20\\
210	21\\
209	22\\
208	23\\
207	24\\
206	25\\
205	26\\
204	27\\
203	28\\
202	29\\
201	30\\
200	31\\
200	32\\
199	33\\
199	34\\
198	35\\
197	36\\
197	37\\
196	38\\
195	39\\
195	40\\
194	41\\
193	42\\
193	43\\
192	44\\
191	45\\
190	46\\
190	47\\
189	48\\
188	49\\
188	50\\
187	51\\
186	52\\
186	53\\
186	54\\
185	55\\
185	56\\
185	57\\
185	58\\
185	59\\
185	60\\
185	61\\
185	62\\
185	63\\
185	64\\
186	65\\
186	66\\
187	67\\
187	68\\
187	69\\
188	70\\
188	71\\
188	72\\
189	73\\
189	74\\
189	75\\
190	76\\
190	77\\
190	78\\
191	79\\
191	80\\
191	81\\
191	82\\
191	83\\
191	84\\
191	85\\
191	86\\
191	87\\
191	88\\
190	89\\
189	90\\
188	90\\
187	90\\
188	91\\
189	92\\
190	93\\
190	94\\
190	95\\
190	96\\
190	97\\
191	98\\
191	99\\
191	100\\
190	101\\
190	102\\
190	103\\
190	104\\
190	105\\
190	106\\
190	107\\
190	108\\
190	109\\
190	110\\
190	111\\
190	112\\
190	113\\
190	114\\
190	115\\
190	116\\
189	117\\
188	117\\
187	116\\
186	115\\
185	115\\
184	114\\
183	113\\
182	112\\
181	112\\
180	111\\
179	110\\
178	109\\
177	109\\
176	108\\
175	107\\
174	107\\
173	106\\
172	105\\
171	105\\
170	104\\
169	103\\
168	103\\
167	103\\
166	102\\
165	102\\
164	101\\
163	101\\
162	100\\
161	100\\
160	99\\
159	99\\
158	98\\
157	98\\
156	98\\
155	97\\
154	97\\
153	96\\
152	96\\
151	95\\
150	95\\
149	94\\
148	94\\
147	93\\
146	93\\
145	93\\
144	92\\
143	92\\
142	91\\
141	91\\
140	90\\
139	90\\
138	89\\
137	89\\
136	88\\
135	88\\
134	87\\
133	87\\
132	86\\
131	85\\
130	85\\
129	84\\
128	84\\
127	83\\
126	82\\
125	82\\
124	81\\
123	80\\
122	80\\
121	79\\
120	78\\
119	77\\
118	77\\
117	76\\
116	75\\
115	74\\
114	73\\
113	72\\
112	71\\
111	70\\
110	69\\
110	68\\
109	67\\
108	66\\
108	65\\
107	64\\
106	63\\
105	62\\
105	61\\
104	60\\
103	59\\
103	58\\
102	57\\
102	56\\
102	55\\
101	54\\
101	53\\
101	52\\
100	51\\
100	50\\
100	49\\
99	48\\
99	47\\
99	46\\
99	45\\
99	44\\
99	43\\
99	42\\
99	41\\
99	40\\
99	39\\
99	38\\
99	37\\
99	36\\
100	35\\
100	34\\
100	33\\
101	32\\
101	31\\
102	30\\
102	29\\
102	28\\
103	27\\
103	26\\
103	25\\
104	24\\
105	23\\
106	22\\
107	21\\
108	20\\
109	19\\
110	18\\
111	17\\
111	16\\
112	15\\
112	14\\
113	13\\
113	12\\
113	11\\
113	10\\
113	9\\
112	8\\
111	8\\
110	8\\
109	7\\
108	7\\
107	7\\
106	7\\
105	7\\
104	7\\
103	7\\
102	7\\
101	7\\
100	7\\
99	7\\
98	7\\
97	7\\
96	8\\
95	8\\
94	8\\
93	8\\
92	8\\
91	9\\
90	9\\
89	10\\
88	10\\
87	11\\
86	12\\
85	12\\
84	13\\
83	14\\
83	15\\
82	16\\
82	17\\
82	18\\
82	19\\
82	20\\
82	21\\
81	22\\
81	23\\
80	24\\
79	24\\
79	25\\
80	26\\
80	27\\
80	28\\
80	29\\
80	30\\
80	31\\
80	32\\
80	33\\
80	34\\
80	35\\
80	36\\
80	37\\
80	38\\
80	39\\
80	40\\
80	41\\
80	42\\
80	43\\
80	44\\
80	45\\
80	46\\
80	47\\
80	48\\
80	49\\
80	50\\
80	51\\
80	52\\
80	53\\
80	54\\
80	55\\
80	56\\
80	57\\
80	58\\
80	59\\
80	60\\
81	61\\
81	62\\
81	63\\
81	64\\
81	65\\
81	66\\
81	67\\
81	68\\
81	69\\
81	70\\
81	71\\
81	72\\
81	73\\
81	74\\
81	75\\
81	76\\
81	77\\
82	78\\
82	79\\
82	80\\
82	81\\
82	82\\
83	83\\
83	84\\
84	84\\
85	85\\
86	86\\
87	86\\
88	87\\
89	88\\
90	88\\
91	89\\
92	90\\
93	91\\
92	92\\
91	92\\
90	92\\
90	93\\
91	93\\
92	94\\
93	95\\
94	96\\
95	97\\
96	98\\
97	99\\
98	100\\
98	101\\
99	102\\
100	103\\
101	104\\
102	105\\
103	106\\
104	107\\
105	108\\
106	109\\
107	110\\
108	111\\
108	112\\
109	113\\
110	113\\
111	114\\
111	115\\
112	116\\
113	117\\
114	118\\
115	119\\
116	120\\
117	121\\
118	122\\
119	123\\
120	124\\
120	125\\
121	126\\
121	127\\
121	128\\
121	129\\
121	130\\
120	131\\
119	130\\
120	131\\
121	132\\
122	133\\
123	134\\
124	135\\
125	136\\
126	137\\
127	138\\
128	139\\
129	140\\
130	141\\
131	142\\
132	143\\
133	144\\
133	145\\
134	146\\
135	147\\
136	148\\
137	149\\
137	150\\
136	151\\
135	150\\
134	150\\
133	150\\
133	151\\
134	152\\
135	153\\
136	154\\
137	155\\
137	156\\
137	157\\
136	158\\
136	159\\
136	160\\
137	161\\
137	162\\
138	163\\
138	164\\
139	165\\
139	166\\
139	167\\
139	168\\
139	169\\
139	170\\
139	171\\
139	172\\
140	173\\
140	174\\
140	175\\
140	176\\
141	177\\
141	178\\
142	179\\
142	180\\
143	181\\
143	182\\
144	183\\
143	184\\
142	184\\
141	183\\
140	182\\
139	182\\
139	183\\
140	184\\
140	185\\
140	186\\
140	187\\
140	188\\
140	189\\
141	190\\
141	191\\
141	192\\
141	193\\
141	194\\
142	195\\
142	196\\
142	197\\
143	198\\
143	199\\
143	200\\
144	201\\
144	202\\
145	203\\
145	204\\
146	205\\
146	206\\
147	207\\
147	208\\
147	209\\
148	210\\
148	211\\
149	212\\
149	213\\
149	214\\
149	215\\
149	216\\
150	217\\
150	218\\
150	219\\
150	220\\
150	221\\
150	222\\
150	223\\
150	224\\
150	225\\
150	226\\
150	227\\
150	228\\
150	229\\
149	230\\
148	231\\
148	232\\
147	233\\
146	234\\
146	235\\
145	236\\
144	237\\
143	237\\
142	238\\
141	239\\
140	240\\
139	241\\
138	242\\
137	243\\
136	244\\
135	245\\
134	246\\
133	246\\
132	247\\
131	248\\
130	249\\
129	250\\
128	251\\
127	252\\
126	253\\
125	254\\
124	255\\
123	256\\
122	256\\
121	257\\
120	258\\
119	259\\
118	260\\
117	261\\
116	262\\
115	263\\
114	264\\
113	265\\
113	266\\
112	267\\
111	268\\
110	269\\
109	270\\
108	271\\
107	272\\
106	273\\
105	274\\
104	275\\
104	276\\
103	277\\
102	278\\
102	279\\
101	280\\
100	281\\
100	282\\
99	283\\
98	284\\
98	285\\
97	286\\
96	287\\
96	288\\
95	289\\
94	290\\
94	291\\
93	292\\
93	293\\
92	294\\
92	295\\
91	296\\
91	297\\
90	298\\
90	299\\
89	300\\
89	301\\
88	302\\
87	303\\
87	304\\
86	305\\
85	306\\
84	307\\
83	308\\
83	309\\
82	310\\
81	311\\
81	312\\
80	313\\
80	314\\
79	315\\
78	316\\
78	317\\
77	318\\
76	319\\
76	320\\
75	321\\
74	322\\
75	322\\
76	321\\
77	322\\
77	323\\
76	324\\
76	325\\
75	326\\
74	327\\
74	328\\
73	329\\
73	330\\
72	331\\
72	332\\
71	333\\
71	334\\
70	335\\
69	336\\
69	337\\
69	338\\
68	339\\
68	340\\
67	341\\
67	342\\
66	343\\
66	344\\
66	345\\
65	346\\
65	347\\
64	348\\
63	349\\
63	350\\
62	351\\
61	352\\
60	353\\
60	354\\
59	355\\
58	356\\
58	357\\
57	358\\
56	359\\
55	359\\
54	360\\
53	361\\
52	362\\
51	363\\
50	364\\
49	365\\
48	366\\
47	366\\
46	367\\
45	368\\
44	368\\
43	369\\
42	369\\
41	370\\
40	371\\
39	371\\
38	372\\
37	373\\
36	373\\
35	374\\
34	375\\
33	375\\
33	376\\
32	377\\
31	378\\
31	379\\
31	380\\
31	381\\
31	382\\
31	383\\
31	384\\
};
\end{axis}
\end{tikzpicture}%